\documentclass[a4paper,10pt]{amsart}
\usepackage{amsmath}
\usepackage[T1]{fontenc}
\usepackage{amssymb}
\usepackage{amsthm}
\newcommand{\md}{{\rm d}}

\newtheorem{thm}{Theorem}[section]
\newtheorem{prop}[thm]{Proposition}
\newtheorem{lem}[thm]{Lemma}
\newtheorem{corol}[thm]{Corollary}

\newtheorem{rem}[thm]{Remark}
\numberwithin{equation}{section}

\title[Jump conditions for Josephson junctions]
{Limiting jump conditions for  Josephson  junctions
in Ginzburg-Landau theory}

\author[Ayman Kachmar]{Ayman Kachmar$^*$}
\thanks{$^*$ Universit\'e Paris-Sud, D\'epartement de math\'ematique,
  B\^at.~425, 91405 Orsay France. {\it E-mail:} {\tt 
ayman.kachmar@math.u-psud.fr}}

\date{\today}

\begin{document}
\maketitle

\begin{abstract}
We consider a  S-N-S Josephson junction modeled through the
Ginzburg-Landau theory. When  the normal material
is sufficiently thin and the applied magnetic
field is below the critical
field of vortex nucleation, we prove to leading order that jump
boundary conditions of the type predicted by de\,Gennes are satisfied
across the junction.
\end{abstract}

\section*{Introduction}
The superconducting proximity effect in a normal metal adjacent to a
superconductor has received a lot of attention by the physics
community, see \cite{Pa} for a review of this phenomenon. This is also the
setting of the Josephson tunneling effect
for superconducting-normal-superconducting junctions (SNS), where a
supercurrent  flows through the normal layer
provided that it is sufficiently thin.\\
The physics literature contains several  approaches to model the
Josephson effect in the frame work of the Ginzburg-Landau theory
of superconductivity. The first modeling in this context is perhaps
due to the physicist
de\,Gennes \cite{deGe}. In the setting of
\cite{deGe}, the complex-valued wave function (whose modulus measures
the density of superconducting electrons) and its derivative are
related linearly on both sides of the normal material, in such a
manner that the supercurrent is conserved through the junction.\\
In this paper, we use a generalized Ginzburg-Landau energy functional
presented in \cite{Chetal}, which
has proved to account rigorously  to various physical
aspects (c.f. \cite{kach3, kach2}). By working in  the London singular limit
(high $\kappa$-regime), we justify asymptotically  the modeling
of \cite{deGe} provided that the applied magnetic field is below the
critical field of vortex nucleation, see
Theorems~\ref{mainthm} \& \ref{mainthm2}.\\
We also mention in this
direction that another justification of the de\,Gennes modeling is
present in a paper of Rubinstein-Schatzman-Sternberg \cite{RuScSt},
who deal with
geometric junctions (weak links) in the framework of the
Ginzburg-Landau theory.\\
We hope to carry out in a forthcoming work   a deeper analysis valid
for higher applied magnetic fields and which provides more details
concerning the supercurent flow and the distribution of vortices in
the junction.
\section{Main results}
We move now to the mathematical set-up of the problem.
Let $\Omega=D(0,1)$ denotes the unit disc in $\mathbb R^2$.
Given $R\in]0,1[$ and $\ell\in]0,R[$,
we introduce the following partition of $\Omega$,
$$\Omega=S\cup \overline N,$$
where
\begin{eqnarray}
&&N=\{x\in\Omega~:~{\rm dist}(x,\partial D(0,R))<\ell\},\label{V-N}\\
&&S_1=D(0,R-\ell),\quad S_2=D(0,1)\setminus D(0,R+\ell),\quad
  S=S_1\cup S_2.\label{V-S}
\end{eqnarray}
We shall suppose that $S$ is the cross section of a cylindrical
superconductor with infinite height and that $N$ is that of a
normal material. By this way, we get a S-N-S Josephson
junction.\\
In Ginzburg-Landau
theory~\cite{GL}, the superconducting properties are described by a
complex valued wave function $\psi$, called the `order parameter',
whose modulus $|\psi|^2$ measures the density of the superconducting
electron Cooper pairs (hence $\psi\equiv0$ corresponds to a normal
state), and a real vector field $ A=(A_1,A_2)$, called the `magnetic
potential', such that the induced magnetic field in the
sample corresponds to ${\rm curl}\, A$. Since,
the superconducting
Cooper electron pairs can diffuse from the superconducting to the
normal material in a normal-superconducting junction, we then have
to consider pairs $(\psi,A)$ defined on
$\Omega$.\\
The basic postulate in the Ginzburg-Landau theory is that the pair
$(\psi,A)$ minimizes the Gibbs free energy, which has
in our case  the following
dimensionless form~\cite{Chetal}~:
\begin{eqnarray}\label{V-EGL}
\mathcal G_{\varepsilon,H}(\psi,A)&=&\int_{\Omega}
|(\nabla-iA)\psi|^2\,\md
x+\frac1{2\varepsilon^2}\int_S(1-|\psi|^2)^2\,\md x\\
&&+\frac{a}{\varepsilon^2}\int_N|\psi|^2\,\md x+\int_\Omega|{\rm
curl}\,A-H|^2\,\md x.\nonumber
\end{eqnarray}
Here,
$\frac1\varepsilon=\kappa$ is a
characteristic of the superconducting material (filling $S$),
$H>0$ is the intensity of the applied magnetic field and
$a>0$ is related to the critical temperature of the material in
$N$. The positive sign of $a$ means that we are above the
critical temperature of the material filling $N$.\\
Minimization of  the
functional (\ref{V-EGL}) will take place in  the space
$$\mathcal H=H^1(\Omega;\mathbb C)\times
H^1(\Omega;\mathbb R^2).$$
We will be interested in the analysis of the asymptotic behaviour of
the minimizers of (\ref{V-EGL}) as $\varepsilon\to0$ (London Limit)
and when the thickness of the ring $N$ is small by taking $
\ell=\ell(\varepsilon)\ll1$ as $\varepsilon\to0$.\\
According to \cite{Gi}, the functional (\ref{V-EGL}) admits a
minimizer $(\psi,A)$ in the space $\mathcal H$.
Our main result is the leading order asymptotic
expansion of the jump of\break $(\nabla-iA)\psi$ across the junction,
i.e. across the boundary of $N$, provided that the order parameter
$\psi$ is not possessing vortices.\\
In order to fix ideas, given a function $f\in H^1(\Omega;\mathbb C)$, we
introduce the {\it jump of} $f$ across $N$ by
\begin{equation}\label{jump}
[f]_N(\theta)= f\left((R+\ell)\,e^ {i\theta}\right)
-f\left((R-\ell)\,e^{i\theta}\right),\quad
\forall~\theta\in[0,2\pi[\,.
\end{equation}

Our  first result concerns the case of very thin rings, of thickness
comparable with $\varepsilon$.

\begin{thm}\label{mainthm}
Let $\Omega=D(0,1)$, $S$ and $N$ as in (\ref{V-N}), (\ref{V-S}), and
$(\psi_\varepsilon,A_\varepsilon)$
be a minimizer of (\ref{V-EGL}). Given $d>0$ and $a>0$,
there exists $\lambda>0$ such that, if the applied
magnetic field satisfies
\begin{equation}\label{Hyp-H}
H\leq \lambda|\ln\varepsilon|\,,\end{equation}
and if $\ell=d\,\varepsilon$, then $|\psi_\varepsilon|>0$ and we have
\begin{equation}\label{jump-cond1}
\lim_{\varepsilon\to0}
\left\|\varepsilon\left[\frac{
n(x)\cdot(\nabla-iA_\varepsilon)\psi_\varepsilon}{\psi_\varepsilon}\right
]_N-
2\sqrt{a}\,\frac{\exp(2\sqrt{a}\,d)-1}{\exp(2\sqrt{a}\,d)+1}
\right\|_{L^2(\mathbb S^1)}=0\,,
\end{equation}
\begin{equation}\label{jump-cond2}
\lim_{\varepsilon\to0}
\big{\|}\,[\psi_\varepsilon
]_N\big{\|}_{L^2(\mathbb S^1)}=0\,.
\end{equation}
Here, $n(x)=\displaystyle\frac{x}{|x|}$ for all $x\in\mathbb
 R^2\setminus\{0\}$, is the unit outward normal
vector of any disc in $\mathbb R^2$.
\end{thm}

\begin{rem}\label{rem-H}
The regime concerning the applied magnetic field $H$ in
Theorem~\ref{mainthm}  corresponds to that below the first critical
field~: When (\ref{Hyp-H}) is satisfied, the order parameter
$\psi_\varepsilon$ has no vortices in $\Omega$.\\
On the other hand, it is well known (c.f. \cite{SaSe})
that there exists $\lambda'>0$ such
that if $H\geq \lambda'|\ln\varepsilon|$,
 the order parameter $\psi_\varepsilon$ has vortices. However, we
are not able to calculate the critical value $\lambda_c$ for which
$$H_{C_1}\sim \lambda_c|\ln\varepsilon|\quad{\rm as}\quad
\varepsilon\to0 \,.$$
This is due to the technical difficulty arising from the very rapid
oscillations of the maximal superconducting density in $N$, see
Section~2 for more details concerning this point.
\end{rem}

\begin{rem}\label{d}
Notice that if one formally makes $d\to0$ in (\ref{jump-cond1}), one
would obtain that the jump across $N$ tends to $0$. This agrees with
experimental and theoretical predictions that the Josephson effect would
be absent in junctions made up of very
thin normal materials, see \cite{Chetal}.
\end{rem}

\begin{rem}\label{deGe-cond}
Let us introduce the vectors
$$X_{\varepsilon}^-=\left(
\begin{array}{c}
\varepsilon\,n(x)\cdot(\nabla-iA_\varepsilon)\psi_\varepsilon
\left(\,(R-\ell)e^{i\theta}\right)\\
\psi_\varepsilon\left(\,(R-\ell)e^{i\theta}\right)
\end{array}\right)\,,$$
$$
X_{\varepsilon}^+=\left(
\begin{array}{c}
\varepsilon\,n(x)\cdot(\nabla-iA_\varepsilon)\psi_\varepsilon
\left(\,(R+\ell)e^{i\theta}\right)\\
\psi_\varepsilon\left(\,(R+\ell)e^{i\theta}\right)
\end{array}\right),\quad\forall~\theta\in[0,2\pi[\,,
$$
and the matrix
$$M_{a,d}=\left(
\begin{array}{cc}
1&
2\sqrt{a}\,\displaystyle\frac{\exp(2\sqrt{a}\,d)-1}{\exp(2\sqrt{a}\,d)+1}\\
0&1
\end{array}\right)\,.$$
Then, (\ref{jump-cond1}) and (\ref{jump-cond2}) can be written in the
equivalent form
\begin{equation}\label{jump-cond3}
\lim_{\varepsilon\to0}
\left\|X_\varepsilon^+-M_{a,d}X_\varepsilon^-\right\|_{L^2(\mathbb
  S^1)}
=0,
\end{equation}
which justifies the boundary condition postulated by de\,Gennes in
\cite{deGe}.
\end{rem}

For thicker rings, we have a result analogous to that of
Theorem~\ref{mainthm}.

\begin{thm}\label{mainthm2}
Let $\Omega=D(0,1)$, $S$ and $N$ as in (\ref{V-N}), (\ref{V-S}), and
$(\psi_\varepsilon,A_\varepsilon)$
be a minimizer of (\ref{V-EGL}). Assume in addition that
$\ell=\ell(\varepsilon)$ satisfies,
\begin{equation}\label{ell}
\varepsilon\ll\ell(\varepsilon)\quad(\varepsilon\to0),\quad\exists\,c>0,~
\forall~\varepsilon\in]0,1],~\ell(\varepsilon)\leq
c\,\varepsilon\ln|\ln\varepsilon|\,.
\end{equation}
Then, given $a>0$,
there exists $\lambda>0$ such that, if the applied
magnetic field satisfies
\begin{equation}\label{Hyp-H}
H\leq \lambda\exp\left(-\frac{2\sqrt{a}\,\ell(\varepsilon)}\varepsilon\right)\,
|\ln\varepsilon|\,,\end{equation}
we have
\begin{equation}\label{jump-cond1'}
\lim_{\varepsilon\to0}
\left\|\varepsilon\left[\frac{
n(x)\cdot(\nabla-iA_\varepsilon)\psi_\varepsilon}{\psi_\varepsilon}\right
]_N-
2\sqrt{a}\,
\right\|_{L^2(\mathbb S^1)}=0\,,
\end{equation}
\begin{equation}\label{jump-cond2'}
\lim_{\varepsilon\to0}
\big{\|}\,[\psi_\varepsilon
]_N\big{\|}_{L^2(\mathbb S^1)}=0\,.
\end{equation}
Here, $n(x)=\displaystyle\frac{x}{|x|}$ for all $x\in\mathbb
 R^2\setminus\{0\}$, is the unit outward normal
vector of any disc in $\mathbb R^2$.
\end{thm}

Notice that the result of Theorem~\ref{mainthm2} agrees with that of
Theorem~\ref{mainthm} when one takes formally $d\to+\infty$ in
(\ref{jump-cond1}).\\

We mention also that the result of Theorem~\ref{mainthm2} is valid
for magnetic fields slightly much lower than that of
Theorem~\ref{mainthm}. Technically, this is due to the fact that we
can not exclude the energy of vortices in $N$ (see Section~2). But
heuristically, the reason is that the maximal superconducting
density (i.e. the positive minimizer of (\ref{V-EGL}) for $H=0$) in
the regime of Theorem~\ref{mainthm2} is small inside the junction,
hence the price of a vortex becomes for magnetic fields less than
that of Theorem~\ref{mainthm}. However, we were not able in this
case to prove that the critical field for vortex nucleation has the
order of $
\exp\left(-\frac{2\sqrt{a}\,\ell(\varepsilon)}\varepsilon\right)\,
|\ln\varepsilon|$, though theoretical predictions in the physics
literature say  that vortices in $N$ would be present for magnetic
fields much below than that of a bulk superconductor, see \cite{Pa}.

The result of Theorem~\ref{mainthm2} is still valid up to lengths
$\ell(\varepsilon)=c_*\,\varepsilon|\ln\varepsilon|$, where
$c_*\in]0,1[$ is  
sufficiently small (this can be checked through minor modifications of
the argument). However, since the corresponding magnetic field will be small
$H\ll1$, we
do not focus on this last regime\footnote{It  is more likely 
that the regime $H\ll1$ is treated, without additional restrictions on
the decay of $H$, by bifurcation arguments. In the present paper,
we treat in detail the
case $H=0$.}.\\

As a by-product of the analysis that we shall carry out, we get a
result concerning the conservation of the current across $N$. The
current is defined as the vector field
\begin{equation}\label{current}
j_\varepsilon=\big{(}i\psi_\varepsilon,(\nabla-iA_\varepsilon)\psi_\varepsilon\big{)}
=\big{(}(i\psi_\varepsilon,(\partial_2-iA_\varepsilon^1)\psi_\varepsilon),
(i\psi_\varepsilon,(\partial_2-iA_\varepsilon^2)\psi_\varepsilon)\big{)}\,,
\end{equation}
where $(\cdot,\cdot)$ denotes the scalar product in $\mathbb C$ when
identified with $\mathbb R^2$.\\
The exact result concerning the current is the following.
\begin{thm}\label{thm-current}
In both regimes of Theorems~\ref{mainthm} \& \ref{mainthm2}, the
circulation of the current is almost conserved across the junction:
\begin{equation}\label{current-conserv}
 \left|\int_{|x|=R+\ell}\tau(x)\cdot j_\varepsilon
-\int_{|x|=R-\ell}\tau(x)\cdot j_\varepsilon\right|=\mathcal
O\left(\varepsilon^{1/2}|\ln\varepsilon|\right)\quad{\rm
as}~\varepsilon\to0\,,
\end{equation}
where $\tau(x)=\frac{x^\bot}{|x|}$ for all $x\in\mathbb R^2\setminus\{0\}$.
\end{thm}

Finally, we comment on some past work concerning non-homogeneous
superconductors and pinning models. Unlike to our situation, pinning
models previously considered correspond to a term in the
Ginzburg-Landau functional having the form $(a-|\psi|^2)^2$, with
$a$ being a smooth function. The first analytic work probably goes
back to \cite{AfSaSe}, where $a=a_\varepsilon$, assumed always
positive, may depend on $\varepsilon$ with the restriction that it
can not oscillate quicker than $|\ln\varepsilon|$ ($|\nabla
a_\varepsilon|\leq C|\ln\varepsilon|$). Compared with our situation,
the discontinuity of the coefficients leads to an order parameter
oscillating between $0$ and $u_\varepsilon$ with $|\nabla
u_\varepsilon|\geq \frac{C}\varepsilon$ in a thin neighborhood of
$\partial S$. Later, in \cite{AnBaPh}, the authors deal with the
case when $a$ has isolated zeros, and prove that vortices appear
first at the zeros of $a$ for magnetic fields having order $1$. In
\cite{AlBr}, the authors allow the function $a$ to have negative
values, but the hypothesis of its smoothness forces the order
parameter to be small on the boundary of the normal side (hence, the
surface superconducting sheath in the normal side is absent).
Moreover, the situation in \cite{AlBr} is more related to the case
of domains with holes and relies on the analysis carried out by the
same authors in \cite{AlBr06}. More recently, the author of the
present paper showed in \cite{kach6} that pinning of vortices is
observed for magnetic fields near the first critical field  
when the function $a$ is a step function.\\
Let us also mention  the very recent work of Alama,
Bronsard \and Sandier in \cite{AlBrSa} where  a
layered superconducting model has been investigated. There, the
expression of the critical field above which vortices are detected
depends strongly on the thickness of the normal regions separating
the superconducting layers.

\subsection*{Organization of the paper}\ \\
Section~2 is devoted to a description of the main points of the
argument.\\
Section~\ref{V-Sec-EnH=0} is devoted to a preliminary analysis of
the minimizers of (\ref{V-EGL}).\\
Section~\ref{V-Sec-meissnerstate} is devoted to an analysis of an
auxiliary variational problem
(this is the variational problem (\ref{Meissner}) describing the Meissner state).\\
Section~\ref{Section-Lowerbound} is devoted to prove a lower bound
of the functional
(\ref{V-EGL}).\\
In Section~\ref{proofs}, we exhibit a vortex-less regime and we
achieve the proofs of Theorems~\ref{mainthm}, \ref{mainthm2} and \ref{thm-current}.\\
Finally, in Appendix~\ref{Appendix-Canon}, we prove a uniqueness
result concerning the solution of the canonical equation without
magnetic field (this is Eq. (\ref{CanonEq})), and in
Appendix~\ref{App-B}, we discuss the difficulty behind the estimate
of the energy of a configuration with vortices on the circle
$\mathbb S_R^1$.

\section{Sketch of  proof}
\subsection*{ Canonical Equation in $\mathbb R^2$
for the case without magnetic field}\ \\
A general technique to tackle asymptotic problems, (successfully
  applied by Helffer \&
  co-authors
for linear problems (cf. e.g. \cite{HeMo3}), and Lu-Pan for non-linear
  problems (cf. e.g. \cite{LuPa96})), is to
  understand the canonical situation. So we consider the case without
  magnetic field, $H=0$, and work with the particular domains,
$$N=\mathbb R\times\,]-\ell,\ell[,\quad S=\mathbb R^2\setminus N\,.$$
When $\ell(\varepsilon)=d\,\varepsilon$, this leads to the equations:
\begin{equation}\label{CanonEq}
\left\{
\begin{array}{l}
-\Delta u=(1-u^2)u\quad {\rm in}~\mathbb R\times\{\mathbb R\setminus[-d,d]\},\\
-\Delta u+a\,u=\quad{\rm in}~\mathbb R\times\,]-d,d[\,,\\
\partial_{x_2} u(\cdot,t_\pm)=\partial_{x_2} u(\cdot, t_\mp),\quad u(\cdot,
 t_\pm)=u(\cdot,t_\mp),\quad t=\pm d\,,
\end{array}
\right.
\end{equation}
Since Eq. (\ref{CanonEq}) arises as a limiting equation of a critical
point of (G-L), we
focus on solutions in the class
$$\mathcal C=\{u\in H^2_{\rm loc}(\mathbb R^2)\cap L^\infty(\mathbb
R^2)~:
~u\geq0\}\,.$$
Using the argument  of \cite[Section~4]{kach3}, Eq. (\ref{CanonEq}) admits in
$\mathcal C$ a unique solution
$$\mathbb R^2\ni(x_1,x_2)\mapsto U(x_2)$$
given by
\begin{equation}\label{CanonEq1}
U(x_2)=\left\{
\begin{array}{l}
\displaystyle\frac{\beta(a,d)\exp(\sqrt{2}\,|x_2|)-1}
{\beta(a,d)\exp(\sqrt{2}\,|x_2|)+1}\quad{\rm if}~ |x_2|\geq d,\\
A(a,d)\left[\exp(\sqrt{a}\,|x_2|)+\exp(-\sqrt{a}\,|x_2|)\right]\quad
{\rm if}~|x_2|<d.\end{array}\right.\end{equation}
where the constants $\beta(a,d)>\exp(-\sqrt{2}\,d)$ and $A(a,d)>0$ are determined explicitly,
but the only important remark is that
$$
\frac{U'(d)}{U(d)}=-\frac{U'(-d)}{U(-d)}=
2\sqrt{a}\frac{\exp(2\sqrt{a}\,d)-1}
{\exp(2\sqrt{a}\,d)+1}\,,
$$
hence the conditions of de\,Gennes are verified
\begin{equation}\label{CanonEq2}
\left(\begin{array}{c}
U'(t)\\
U(t)\end{array}\right)_{t=d}=\left(
\begin{array}{cc}
1&2\sqrt{a}\,\displaystyle\frac{\exp(2\sqrt{a}\,d)-1}
{\exp(2\sqrt{a}\,d)+1}\\
0&1\end{array}\right)\,
\left(\begin{array}{c}
U'(t)\\
U(t)\end{array}\right)_{t=-d}\,.
\end{equation}
\subsection*{The case without magnetic field in a bounded domain}\ \\
Now we return to minimizing (\ref{V-EGL}) when $H=0$. We prove in
\cite[Theorem~1.1]{kach3} that (\ref{V-EGL}) has, up to a gauge
transformation, a unique minimizer $(u_\varepsilon,0)$ where
$u_\varepsilon\in H^2(\Omega)$ is a real-valued function,
$0<u_\varepsilon<1$ in $\overline\Omega$ (for $\varepsilon$ small
enough), and $u_\varepsilon$ solves the equation
$$-\Delta
u_\varepsilon=\frac{1}{\varepsilon^2}(1-u_\varepsilon^2)u_\varepsilon\,
1_S-\frac{a}{\varepsilon^2}u_\varepsilon\, 1_N\,.$$
Then, by a  blow-up argument, we generalize (\ref{CanonEq2})
asymptotically as $\varepsilon\to0$, see Theorem~\ref{mainthm-H=0}.
\subsection*{The case with magnetic field: A vortex-less regime}\ \\
Now we return to minimizers $(\psi,A)$ of (\ref{V-EGL}) when
$H>0$. Following an idea of Lassoued-Mironescu \cite{LaMi}, we
introduce a normalized density
$$\varphi=\frac{\psi}{u_\varepsilon}\,.$$
Then, $|\varphi|\leq1$ and we are led to work with the functional
(see Theorem~\ref{V-lem-psi<u}):
\begin{eqnarray*}
&&\hskip-0.5cm\mathcal F_{\varepsilon,H}(\varphi,A)=\int_\Omega
u_\varepsilon^2
|(\nabla-iA)\varphi|^2\,\md x\\
&&\hskip2.5cm+\frac1{2\varepsilon^2} \int_S
u_\varepsilon^4(1-|\varphi|^2)^2\,\md x+\int_\Omega |{\rm
curl}\,A-H|^2\,\md x.\nonumber
\end{eqnarray*}
Then, following \cite{Sand, SaSe00}, we construct a family of
disjoint balls $(B(a_i,r_i))_i$ (see
Proposition~\ref{V-lem-vortexballs}) such that
$\sum_ir_i\leq |\ln\varepsilon|^{-10}$ and
$$\{x\in\overline\Omega~:~|\varphi(x)|\leq
    1-|\ln\varepsilon|^{-2}\}\subset
\bigcup_i B(a_i,r_i)\,.$$
This permits to obtain, for a given number $\alpha\in]0,\frac12[$, 
a lower
    bound of the energy (see Theorem~\ref{LowerBound-thm}):
\begin{equation}\label{Lowerbound-Int}
\mathcal F_{\varepsilon,H}(\varphi,A)\geq
M_0(\varepsilon,H)+\sum_i \left(2\pi\alpha\,
m_i(\varepsilon)|\ln\varepsilon|-2H\right)d_i-C\,H|\ln\varepsilon|^{-10}\,,
\end{equation}
where $C>0$ is an explicit constant, $d_i$ is the degree of
$\varphi/|\varphi|$ on $\partial B(a_i,r_i)$,
$$m_i(\varepsilon)=\min_{x\in B(a_i,r_i)}u_\varepsilon^2(x)\,,$$
and
\begin{equation}\label{Meissner}
M_0(\varepsilon, H)=\inf_{A\in H^1(\Omega;\mathbb R^2)}\mathcal
F_{\varepsilon,H}(1,A)\,.
\end{equation}
By comparing (\ref{Lowerbound-Int}) with  what is existing in the
literature (c.f. \cite{SaSe}), we suspect that the lower bound
(\ref{Lowerbound-Int}) is not optimal in the sense that $\alpha$
should be equal to $1$. This restriction is actually due to the
particular expression of $\mathcal F_{\varepsilon,H}$,
where a penalization term for $|\varphi|$
is absent in the energy of $N$.\\
The infimum  in (\ref{Meissner}) is achieved by a unique vector field
$\mathcal A=\frac{H}{u_\varepsilon^2}\nabla^\bot h_\varepsilon$, where
$h_\varepsilon ~:\Omega\mapsto]0,1[$ satisfies a London equation with
    weight (see (\ref{V-hepsilon'})). Thus,  we get the upper
bound
$$\mathcal F_{\varepsilon,H}(\varphi,A)\leq M_0(\varepsilon,H).$$
When this upper bound is matched with (\ref{Lowerbound-Int}), we
deduce that  all the $d_i's$ are equal to $0$ provided that
$$H\leq \alpha\pi \left(\inf_i
m_i(\varepsilon)\right)|\ln\varepsilon|\,.$$
If $\ell(\varepsilon)\gg \varepsilon$ and $|a_i|=R$, we have by
Theorem~\ref{L-infty-estimate} that
$$m_i(\varepsilon)\approx
\exp\left(-\frac{2\sqrt{a}\,\ell(\varepsilon)}{\varepsilon}\right)\,.$$
Since we can not exclude the possibility of a vortex ball centered
on the circle $\mathbb S_R^1$, we restrict ourselves when
$\ell(\varepsilon)\gg\varepsilon$  to magnetic fields $H$ satisfying
$$H\leq \lambda
\exp\left(-\frac{2\sqrt{a}\,\ell(\varepsilon)}{\varepsilon}\right)|\ln\varepsilon|$$
in order to insure the absence of vortices.\\
Now, in the absence of vortices we get an energy estimate (see Theorem~\ref{existence-lambda})
$$\int_\Omega\left(|(\nabla-iA)\varphi|^2
+|{\rm curl}\,A-H\,h_\varepsilon|^2\right)\md x+\frac1{\varepsilon^2}
\int_S(1-|\varphi|^2)^2\md x\ll1\quad(\varepsilon\to0)\,.$$
Then, we implement $L^2$-estimates for the equations of $\varphi$
in order to deduce that
$$[\varphi]_N\to 0\quad{\rm in}~L^2(\mathbb S^1),\quad |\varphi|\to 1,\quad
n(x)\cdot(\nabla-iA)\varphi\to 0\quad{\rm in}~L^2(\partial N)\,,$$
which permits us to  deduce
Theorems~\ref{mainthm} and \ref{mainthm2}, see
Section~\ref{existence-lambda}.

\section{Preliminary analysis of minimizers}
\label{V-Sec-EnH=0}

\subsection{The case without applied magnetic field}
This section is devoted to a summary of the main results obtained in
\cite{kach3} which deal with minimizers of (\ref{V-EGL}) when the
applied magnetic field $H=0$.\\
We keep the notation introduced in Section~1. Upon taking $A=0$ and
$H=0$ in (\ref{V-EGL}), one is led to introduce the functional
\begin{equation}\label{V-EnH=0}
\mathcal G_\varepsilon(u):=\int_\Omega|\nabla u|^2\,\md x
+\frac{1}{2\varepsilon^2} \int_S(1-u^2)^2\,\md
x+\frac{a}{\varepsilon^2} \int_Nu^2\,\md x,
\end{equation}
defined for functions in $H^1(\Omega;\mathbb R)$.\\
We introduce
\begin{equation}\label{V-C0}
C_0(\varepsilon)=\inf_{u\in H^1(\Omega;\mathbb R)}
\mathcal G_\varepsilon(u)\,.
\end{equation}

The next theorem is a summary of Theorem~1.1 in \cite{kach3}.

\begin{thm}\label{V-thm-kach3}
Given $a>0$ and $d>0$, there exists $\varepsilon_0$ such that for all
$\varepsilon\in]0,\varepsilon_0[$, the functional (\ref{V-EnH=0})
admits in $H^1(\Omega;\mathbb R)$ a
minimizer $u_\varepsilon\in C^2(\overline{S})\cup C^2(\overline N)$
such that
$$0<u_\varepsilon<1\quad {\rm in}~\overline{\Omega}.$$
Furthermore, with our choice of the domains $\Omega,N$ and $S$ in
(\ref{V-N}) and (\ref{V-S}), the function $u_\varepsilon$ is radial.\\
If $H=0$, minimizers of (\ref{V-EGL}) are gauge equivalent to the
state $(u_\varepsilon,0)$.
\end{thm}

Let us just mention why we focus only on the regime $\varepsilon\to 0$ 
in the statement of  Theorem~\ref{V-thm-kach3}. Notice that $u\equiv
0$ is a critical point of the functional (\ref{V-EnH=0}), so we would
like to exclude the possibility that this critical point is stable. To
this end, 
we define the following eigenvalue~:
\begin{eqnarray*}
&&\lambda_1(a,d,\varepsilon)=\inf\left\{\int_{S}\left(|\nabla\phi|^2
-\frac1{\varepsilon^2}|\phi|^2\right)\md x\right.\\
&&\left.+\int_{N}\left(
|\nabla\phi|^2+\frac{a}{\varepsilon^2}|\phi|^2\right)\md
x~:\quad \phi\in
H^1(\Omega),\,\|\phi\|_{L^2(\Omega)}=1\right\}.\nonumber
\end{eqnarray*}
Then when $\lambda_1(a,d,\varepsilon)<0$, the corresponding
eigenfunction of $\lambda_1(a,d,\varepsilon)$ provides us with a test
configuration whose energy is below that of $u\equiv0$. On the other
hand, this last condition of the sign of $\lambda_1(a,m,\varepsilon)$
is easily verified when $\varepsilon\to0$, thanks in particular to the
min-max principle.\\ 

Let us recall the notation of the {\it jump} across $N$ introduced in
(\ref{jump}). Using a blow-up argument and a result concerning
uniqueness of critical points of the
functional (\ref{V-EnH=0}) in the model case of the entire plane,
we establish now Theorem~\ref{mainthm}
in the case when $H=0$.

\begin{thm}\label{mainthm-H=0}
Let $u_\varepsilon$ be the positive minimizer of (\ref{V-EnH=0})
introduced in Theorem~\ref{V-thm-kach3}. Then, if
$\ell(\varepsilon)=d\,\varepsilon$, we have
\begin{equation}\label{jump-cond1-H=0}
\lim_{\varepsilon\to0}
\left\|\varepsilon\left[\frac{
n(x)\cdot\nabla u_\varepsilon}{u_\varepsilon}\right
]_N-
2\sqrt{a}\,\frac{\exp(2\sqrt{a}\,d)-1}{\exp(2\sqrt{a}\,d)+1}
\right\|_{L^\infty(\mathbb S^1)}=0\,,
\end{equation}
\begin{equation}\label{jump-cond2-H=0}
\lim_{\varepsilon\to0}
\big{\|}\,[u_\varepsilon
]_N\big{\|}_{L^\infty(\mathbb S^1)}=0\,,
\end{equation}
\begin{equation}\label{Bnd}
\lim_{\varepsilon\to0}\left(\sup_{\theta\in[0,2\pi[}
\left|u_\varepsilon\left(\,(R\pm\ell)e^{i\theta}\right)
-\widetilde A(a,d)\right|\right)=0\,.
\end{equation}
Here, $n(x)=\displaystyle\frac{x}{|x|}$ for all $x\in\mathbb
 R^2\setminus\{0\}$, and $A(a,d)>0$ is an explicit
 constant\footnote{The expression of $\widetilde A(a,d)$ is given explicitly in
 the Appendix.}.\\
On the other hand, if $\ell(\varepsilon)\gg\varepsilon$ (this
covers the regime (\ref{ell})), then
 we have,
\begin{equation}\label{jump-cond1-H=0'}
\lim_{\varepsilon\to0}
\left\|\varepsilon\left[\frac{
n(x)\cdot\nabla u_\varepsilon}{u_\varepsilon}\right
]_N-
2\sqrt{a}\,
\right\|_{L^\infty(\mathbb S^1)}=0\,,
\end{equation}
\begin{equation}\label{Bnd'}
\lim_{\varepsilon\to0}\left(\sup_{\theta\in[0,2\pi[}
\left|u_\varepsilon\left(\,(R\pm\ell)e^{i\theta}\right)
-A(a)\right|\right)=0\,,
\end{equation}
where $A(a)>0$ is an explicit constant.
\end{thm}
\paragraph{\bf Proof.}
Let us treat the case when $\ell(\varepsilon)=d\,\varepsilon$,
$d>0$. Let $(r,\theta)$ be polar coordinates, $0<r<1$,
$-\pi<\theta<\pi$,  and set
$$t=R-r,\quad s=R\,\theta\,.$$
Given $s_0\in[-R\pi,R\pi[$,
we define the rescaled function,
$$\widetilde u_\varepsilon(s,t)=u_\varepsilon\left((R-\varepsilon t)
e^{i\varepsilon
(s-s_0)/R}\right)\,,\quad  \frac{R-1}\varepsilon<t<\frac{1-R}\varepsilon,
~-\pi\frac{R}\varepsilon<s-s_0<\pi\frac{R}\varepsilon.$$
The equation of $\widetilde u_\varepsilon$ becomes:
$$\left\{
\begin{array}{l}
-\Delta_\varepsilon\, \widetilde u_\varepsilon =(1-\widetilde
u_\varepsilon^2)\widetilde u_\varepsilon,
\quad 0<t<\frac{1-R}\varepsilon,~|s-s_0|<\pi\frac{R}\varepsilon,
\\ \\
-\Delta_\varepsilon\,\widetilde u_\varepsilon+am\,\widetilde
u_\varepsilon=0,
\quad \frac{R-1}\varepsilon<t<0,~|s-s_0|<\pi\frac{R}\varepsilon,\\ \\
\displaystyle\frac{\partial\widetilde u_\varepsilon}{\partial
t}(\cdot,t_\pm)=\displaystyle\frac{\partial\widetilde
u_\varepsilon}{\partial t}(\cdot,t_\mp),\quad
\widetilde u_\varepsilon(\cdot,t_\pm)=\widetilde u(\cdot,t_\mp)~
\text{for }t=\pm d\,,
\end{array}\right.$$
where
$$\Delta_\varepsilon=\left(1-\varepsilon
\frac{t}R\right)^{-2}\partial_s^2
+\partial_t^2 -\frac{\varepsilon }{\left(R-\varepsilon t\right)}\partial_t.$$
Now, by elliptic estimates, the function $\widetilde u_\varepsilon$
converges to a function $u$ in $W^{2,\infty}_{\rm loc}(\mathbb
R^2)$. Furthermore, $u$ solves (\ref{CanonEq}) in $\mathcal C$, and
by \cite[Lemma~5.2]{kach3}, there exist constants $k_0,c_0>0$ such
that $u(0,k_0)>c_0$. 
Thus, we conclude by Theorem~\ref{CanonEq-thm} that $u(s,t)=
U(t)$, where $U$ is given in (\ref{CanonEq2}), and therefore, by
coming back to the initial scale,
$$\forall~k\in\{0,1,2\},\quad
\lim_{\varepsilon\to 0}\varepsilon^k\left\|u_\varepsilon(s,t)
-U\left(\frac t\varepsilon\right)\right\|_{W^{k,\infty}(\{|s-s_0|\leq \pi
  R\varepsilon,\,|t|\leq (1-R)\varepsilon\})}=0,$$
and the convergence is uniform with repect to $s_0\in[-\pi R,\pi R[$.
This yields (\ref{jump-cond1-H=0})-(\ref{Bnd}).\\
The statements concerning the case
when $\ell(\varepsilon)\gg\varepsilon$ follows from
\cite[(5.20)]{kach3}.\hfill$\Box$\\

We shall need the following remarkable properties of
$u_\varepsilon$, that distinguish the different regimes considered
in this paper.

\begin{lem}\label{u>0}
With the notations and hypotheses of Theorem~\ref{mainthm-H=0}, if $\ell=d\,\varepsilon$,
there exists an explicit constant $c(a,d)>0$ such that
$$u_\varepsilon(x)>c(a,d),\quad\forall~x\in\overline\Omega\,.$$
\end{lem}

For the case when $\ell(\varepsilon)$ satisfies (\ref{ell}), the
asymptotic behavior of $u_\varepsilon$ becomes completely different in
the sense that it is close to $0$ inside $N$.\\
In order to be precise we introduce the function~:
\begin{equation}\label{V-U(t)}
V(t)=\frac{\beta\exp(\sqrt{2}\,t)-1}{\beta\exp(\sqrt{2}\,t)+1}\quad
(t\geq 0),\quad V(t)=A\,\exp(\sqrt{a}\,t)\quad (t<0),
\end{equation}
together with
 the `signed distance' to
the boundary of $S$,
\begin{equation}\label{V-tS}
t_S(x)= {\rm dist}(x,\partial S)\quad (x\in S),\quad t_S(x)=-{\rm
dist}(x,\partial S)\quad  (x\in D(0,1)\setminus S).
\end{equation}
Here the constants $\beta$ and $A$ are given by~:
\begin{equation}\label{VI-l,bet}
\beta=\frac{\sqrt{2}+\sqrt{a+2}}{\sqrt{a}},\quad
A=\frac{\sqrt{2}+\sqrt{a+2}-\sqrt{a}}{\sqrt{2}+\sqrt{a+2}+\sqrt{a}}.
\end{equation}

\begin{thm}\label{L-infty-estimate}
Assume that $\ell(\varepsilon)$ satisfies (\ref{ell}). Then, we have
\begin{equation}\label{V-u=v'=}
\left\|u_\varepsilon-V\left(\frac{t_S(x)}{\varepsilon}\right)\right\|
_{L^\infty(\Omega)}=o(1)\quad (\varepsilon\to0),
\end{equation}
where the functions $t_S$ and $V$  have been introduced in
(\ref{V-U(t)})-(\ref{V-tS}).\\
Moreover, there exist a positive constant $\varepsilon_0$ and a
function $]0,1]\ni \varepsilon\mapsto g(\varepsilon)\in]0,1[$ such
that
 $g(\varepsilon)\ll1$
and  for all $\varepsilon\in]0,\varepsilon_0]$, one has the estimate
\begin{equation}\label{2sidedEstimate}
(A-g(\varepsilon))\exp\left(\frac{\sqrt{a}\,t_S(x)}\varepsilon\right)\leq
u_\varepsilon(x)\leq (A+g(\varepsilon))\exp\left(\frac{\sqrt{a}\,
  t_S(x)}\varepsilon\right),\quad\forall~x\in \overline N.
\end{equation}
\end{thm}
\paragraph{\bf Proof.}
The asymptotic behavior in (\ref{V-u=v'=}) has been obtained in
\cite{kach3}.
We have only to prove the improved estimate in $N$, i.e.
(\ref{2sidedEstimate}).\\
Let us show how one can obtain the lower bound. Let us introduce the
function~:
$$v_\varepsilon(x)=C\exp\left(\frac{\delta\,
  t_S(x)}\varepsilon\right),$$
where $C$ and $\delta$ are positive constants to be specified later.\\
Let us recall that by definition, the function $t_S$ is written as
$$t_S(x)=\left\{\begin{array}{cl}
R-\ell(\varepsilon)-|x|&{ \rm if}~R-\ell(\varepsilon)\leq |x|\leq R,\\
|x|-R-\ell(\varepsilon)&{\rm if }~R<|x|\leq R+\ell(\varepsilon),
\end{array}\right.$$
where the constant  $R\in]0,1[$ has been introduced for  defining
$S$ and $N$, see (\ref{V-N}) and (\ref{V-S}).\\
Consequently, the function $t_S$ is smooth in each of the following
two annuli of $N$:
$$N_-=\{x\in N~:~R-\ell(\varepsilon)<|x|<R\},\quad N_+=\{x\in N~:~
R<|x|<R+\ell(\varepsilon)\}.$$ One then checks easily that
\begin{eqnarray*}
-\Delta(u_\varepsilon-v_\varepsilon)+\frac{a}{\varepsilon^2}
(u_\varepsilon-v_\varepsilon) &=&
\frac{\delta^2}{\varepsilon^2}\left[1
-\frac{a}{\delta^2}\pm\frac{\varepsilon}{\delta}|x|^{-1}\right]v_\varepsilon\\
&\geq&\frac{\delta^2}{\varepsilon^2}\left[1
-\frac{a}{\delta^2}-\frac{\varepsilon}{\delta(R-\ell(\varepsilon))}
\right]v_\varepsilon \quad {\rm in}~N_\pm.
\end{eqnarray*}
It is  a result of the  asymptotic
formula (\ref{V-u=v'=}) that there exist a constant
$\varepsilon_0$ and  a function
$]0,\varepsilon_0]\ni\varepsilon\mapsto f(\varepsilon)\in]0,1[$ such
that $f(\varepsilon)\ll1$ as $\varepsilon\to0$ and
$$\left|{u_\varepsilon}_{|_{\partial N}}-A\right|\leq f(\varepsilon),
\quad\varepsilon\in]0,\varepsilon_0].$$
Therefore, gathering all the
above remarks, we get for
$$\delta=\sqrt{a+\frac{\varepsilon^2}{4(R-\ell(\varepsilon))^2}
-\frac{\varepsilon} {2(R-\ell(\varepsilon))}}\,,\quad
C=A-2f(\varepsilon)$$ and when $\varepsilon\in]0,\varepsilon_0]$,
\begin{equation}\label{appl-MP}
\left\{
\begin{array}{cl}
-\Delta(u_\varepsilon-v_\varepsilon)
+\displaystyle\frac{a}{\varepsilon^2}(u_\varepsilon-v_\varepsilon)\geq
 0&{\rm in}~N_\pm\\
u_\varepsilon(x)-v_\varepsilon(x)>0&{\rm on}~\partial N\,.
\end{array}
\right.\end{equation}
Two cases may occur regarding the gradient of $u_\varepsilon$ on the
circle $|x|=R$,  either $u'_\varepsilon(R)\leq0$ or
$u'_\varepsilon(R)>0$.\\ 
If $u'_\varepsilon(R)\leq0$, then we get in addition to
(\ref{appl-MP})
$$\displaystyle\frac{\partial}{\partial\nu_{N_+}}(u_\varepsilon-v_\varepsilon)>0\quad{\rm
  on}~(\partial N_+)\cap N.$$
Here, we recall that $\nu_{N_\pm}$ denote the unit outward
normal vectors of
the boundaries of $N_\pm$.\\
Therefore, the application of the strong maximum principle yields
that
$$u_\varepsilon-v_\varepsilon\geq0\quad{\rm  in}\quad
\overline{N_+}\,.$$
This last lower bound when combined with (\ref{appl-MP}) yields
$$ \left\{
\begin{array}{cl}
-\Delta(u_\varepsilon-v_\varepsilon)
+\displaystyle\frac{a}{\varepsilon^2}(u_\varepsilon-v_\varepsilon)\geq
 0&{\rm in}~N_-\\
u_\varepsilon(x)-v_\varepsilon(x)>0&{\rm on}~\partial N_-\,.
\end{array}\right.$$
Hence by the strong maximum principle,
$u_\varepsilon-v_\varepsilon\geq 0$ in $\overline{N_-}$. Therefore, we
deduce that
$$u_\varepsilon-v_\varepsilon\geq 0\quad{\rm in}~\overline{N}\,,$$
which is nothing but the lower bound
we wish to prove for the function $u_\varepsilon$. The same argument
holds when $u'_\varepsilon(R)>0$, but by changing the roles of $N_+$
and $N_-$.\\
The proof of the upper bound follows the same lines above, so we
omit the details.\hfill$\Box$\\

\subsection{The case with magnetic field}\label{V-sec-minimizers}
This section is devoted to a preliminary analysis of the minimizers
of (\ref{V-EGL}) when $H\not=0$. The main point that we shall show
is how to extract the singular term $C_0(\varepsilon)$ (cf.
(\ref{V-C0})) from the energy of a minimizer.

Notice that the existence of minimizers is standard starting from a
minimizing sequence (cf. e.g. \cite{Gi}). A standard choice of gauge
permits one to assume that the magnetic potential satisfies
\begin{equation}\label{V-gauge}
{\rm div}\,A=0\quad {\rm in}~\Omega,\quad \nu\cdot A=0\quad{\rm
on}~\partial\Omega,
\end{equation}
where $\nu$ is the outward unit normal vector of
$\partial\Omega$.\\
With this choice of gauge, one is able to prove (when the boundaries
of $\Omega$ and $N$ are  smooth) that a minimizer $(\psi,A)$ is in $
C^1(\overline\Omega;\mathbb C)\times C^1(\overline\Omega;\mathbb
R^2)$. One  has also the following regularity (cf.
\cite[Appendix~A]{kach3}),
$$\psi \in C^2(\overline S;\mathbb C)\cup C^2(\overline N;\mathbb
C),\quad A\in C^2(\overline S;\mathbb R^2)\cup C^2(\overline
N;\mathbb R^2).$$

The next lemma is inspired from the work of Lassoued-Mironescu (cf.
\cite{LaMi}).

\begin{lem}\label{V-lem-psi<u}
Let $(\psi,A)$ be a minimizer of (\ref{V-EGL}). Then
$0\leq|\psi|\leq u_\varepsilon$ in $\Omega$, where $u_\varepsilon$
is the positive
minimizer of (\ref{V-EnH=0}).\\
Moreover, putting $\varphi=\frac{\psi}{u_\varepsilon}$, then the
energy functional (\ref{V-EGL}) splits in the form~:
\begin{equation}\label{V-splittingEn}
\mathcal G_{\varepsilon,H}(\psi,A)=C_0(\varepsilon)+ \mathcal
F_{\varepsilon,H}(\varphi,A),\end{equation} where $C_0(\varepsilon)$
has been introduced in (\ref{V-C0}) and the new functional $\mathcal
F_{\varepsilon,H}$ is defined by~:
\begin{eqnarray}\label{V-reducedfunctional}
&&\hskip-0.5cm\mathcal F_{\varepsilon,H}(\varphi,A)=\int_\Omega
u_\varepsilon^2
|(\nabla-iA)\varphi|^2\,\md x\\
&&\hskip2.5cm+\frac1{2\varepsilon^2} \int_S
u_\varepsilon^4(1-|\varphi|^2)^2\,\md x+\int_\Omega |{\rm
curl}\,A-H|^2\,\md x.\nonumber
\end{eqnarray}
\end{lem}
\paragraph{\bf Proof}\ \\
The equality (\ref{V-reducedfunctional}) results from a direct but some
how long calculation, which permits to deduce in particular that
$\varphi$ is a solution of the equation
$$-(\nabla-iA)u_\varepsilon^2(\nabla-iA)\varphi=
1_S\frac{u_\varepsilon^4}{\varepsilon^2}(1-|\varphi|^2)^2\varphi\,.$$\\
{\it Proof of $|\psi|\leq u_\varepsilon$.}\\
It is sufficient to prove that $|\varphi|\leq 1$. We shall invoke an
energy argument which we take from \cite{DGP}.\\
Let us introduce the set
$$\Omega_+=\{x\in\overline\Omega~:~|\varphi(x)|>1\}\,,$$
together with the functions (defined in $\Omega_+$)~:
$$f=\frac{\varphi}{|\varphi|}\,,\quad
\widetilde\varphi=[\,|\varphi|-1]_+f\,.$$
Then, it results from a direct calculation together with the
weak-formulation of the equation satisfied by $\varphi$ that
\begin{eqnarray*}
0&=&\int_{\Omega_+}
\left(|\nabla|\varphi|\,|^2+(|\varphi|-1)
|\varphi|\,|(\nabla-iA)f|^2\right)u_\varepsilon^2\,\md x\\
&&+\frac1{\varepsilon^2}\int_{\Omega_+\cap S}
\left(1+|\varphi|)(1-|\varphi|)^2|\varphi|\right)u_\varepsilon^4\,\md
x\,.
\end{eqnarray*}
Therefore, this yields that $\Omega_+\subset N$ and that
$|\nabla|\varphi|\,|\equiv 0$ in $\Omega_+$. Hence, $|\varphi|$ is
constant in each connected component of $\Omega_+$, which shows that
$|\varphi|\equiv1$ in $\Omega_+$. This contradicts the definition of
$\Omega_+$ unless $\Omega_+=\emptyset$.\hfill$\Box$\\

The estimate of the next lemma is very useful for exhibiting a
vortex-less regime for  minimizers  of (\ref{V-EGL}).

\begin{lem}\label{grad-Abd}
Let $(\psi,A)$ be a minimizer of (\ref{V-EGL}). There exist
constants $C>0$ and $\varepsilon_0\in]0,1]$ such that, if the applied magnetic
field satisfies $H\ll\frac1{\varepsilon}$, then we have
$$|(\nabla-iA)\psi|\leq\frac{C}\varepsilon,\quad\forall\varepsilon\in]0,\varepsilon_0]\,.$$
\end{lem}
\paragraph{\bf Proof}
The proof is essentially due to B\'ethuel-Rivi\`ere \cite{BeRi}, but we include
the main steps for the reader's convenience.\\
Since $|(\nabla-iA)\psi|$ is a gauge invariant quantity, we assume
that we
are in the Coulomb gauge (\ref{V-gauge}).\\
Let us assume that the conclusion of the lemma were false. Then there
exists a subsequence, denoted again $\{\varepsilon\}$ and points
$(x_\varepsilon)\in\Omega$ such that
\begin{equation}\label{contradiction}
\varepsilon|(\nabla-iA)\psi(x_\varepsilon)|\to\infty\,.
\end{equation}
We define the rescaled functions
$$v_\varepsilon(x)=\psi(x_\varepsilon+\varepsilon x),\quad
B_\varepsilon(x)=\varepsilon A(x_\varepsilon+\varepsilon x)\,,$$
together with the rescaled domain
$$\Omega_\varepsilon=(\Omega-x_\varepsilon)/\varepsilon\,.$$
Notice that $(v_\varepsilon,B_\varepsilon)$ satisfies the equations
$$\left\{\begin{array}{rcl}
-\Delta
v_\varepsilon+2iB_\varepsilon\cdot\nabla&=&\left[(1-|v_\varepsilon|^2)\,1_S-a\,1_N+|B_\varepsilon|^2
\right]v_\varepsilon\quad{\rm in}~
\Omega_\varepsilon,\\
-\Delta B_\varepsilon&=&\varepsilon^2\left(iv_\varepsilon\,,\,\nabla
v_\varepsilon-iB_\varepsilon v_\varepsilon\right)\quad
{\rm in}~
\Omega_\varepsilon,\\
{\rm curl}\,B_\varepsilon&=&\varepsilon^2H^2\quad {\rm on}~
\partial\Omega_\varepsilon,\\
n\cdot\nabla v_\varepsilon&=&0\quad
{\rm on}~
\partial\Omega_\varepsilon\,.
\end{array}\right.$$
Notice that $v_\varepsilon$ is a weak solution in $\Omega$ of
the first equation above because of the transmission conditions
$$\mathcal T_{\partial N}^{\rm int}(n\cdot\nabla
v_\varepsilon)=\mathcal T_{\partial N}^{\rm ext}(n\cdot\nabla v_\varepsilon),$$
where
$$\mathcal T_{\partial N}^{\rm int}~:~H^1(N)\mapsto L^2(\partial N),\quad
\mathcal T_{\partial N}^{\rm int}~:~H^1(\Omega\setminus \overline N)
\mapsto L^2(\partial N),$$
are the `interior' and `exterior' trace operators.\\
With the choice of gauge in (\ref{V-gauge}), we get by Poincar\'e's
Lemma that $\|A\|_{H^1(\Omega)}\leq C\|{\rm
curl}\,A\|_{L^2(\Omega)}$. Then, by $L^2$ elliptic estimates and the
assumption $H\ll \frac1{\varepsilon}$, we have
$$\|A_\varepsilon\|_{H^2(\Omega)}\leq C\|{\rm
  curl}\,A\|_{L^2(\Omega)}\ll
\frac1\varepsilon.$$
Hence, by the Sobolev
embedding theorem, we have in  the new scale,
$$\lim_{\varepsilon\to0}\|B_\varepsilon\|_{L^\infty(\Omega_\varepsilon)}=0\,.$$
Now, by Lemma~\ref{V-lem-psi<u}, $|v_\varepsilon|\leq
u_\varepsilon\leq 1$, hence, the right hand side of the equation of
$v_\varepsilon$ becomes bounded. Therefore,
$$\|\Delta v_\varepsilon\|_{L^p(B_R)}\leq C_R+
2\|B_\varepsilon\|_{L^\infty(\Omega_\varepsilon)}\times\|\nabla v_\varepsilon\|_{L^p(B_R)},\quad\forall~p>2\,,$$
where $B_R$ is any fixed ball of radius $R$.\\
By elliptic regularity theory, $v_\varepsilon$ becomes bounded in
$W^{2,p}(B_R)$ for all $p>2$, hence,  by the Sobolev embedding
theorem,  in $C^{1,\alpha}(B_R)$ for any $\alpha\in]0,1[$. Since
$C^{1,\alpha}(B_R)$ is compactly embedded in $C^1(B_R)$, we get by a
diagonal argument the existence of  a function $v\in C^1(\mathbb R
^2)$ such that, upon extraction of a subsequence, $v_\varepsilon$
converges
to $v$ locally in $C^1$.\\
Now, we know from the equation of $B_\varepsilon$ that $\Delta
B_\varepsilon$ is locally bounded in $L^\infty$. So again, by
elliptic estimates, and since $\|B_\varepsilon\|_{L^\infty}\to0$ as
$\varepsilon\to0$, we get upon extraction of a
subsequence that $B_\varepsilon$ converges to $0$ locally in $C^1$.\\
Therefore, we get, by returning to the initial scale,
$$\varepsilon|(\nabla-iA)\psi(x_\varepsilon)|=|(\nabla-iB_\varepsilon)v_\varepsilon(0)|$$
is convergent, hence contradicting
(\ref{contradiction}).\hfill$\Box$\\

Now, Lemma~\ref{grad-Abd} permits to conclude the following result.

\begin{lem}\label{BBH-thm3.3}
Assume that $(\psi,A)$ is a minimizer of (\ref{V-EGL}) and
let $\varphi=\frac\psi{u_\varepsilon}$. There exists a constant
$\mu_0>0$ such that if
$$\frac1{\varepsilon^2}\int_S (1-|\varphi|^2)^2\,\md x\leq \mu_0\,,$$
then $|\varphi|\geq \frac12$ in $S$.
\end{lem}
\paragraph{\bf Proof.}
Lemma~\ref{grad-Abd} and the diamagnetic inequality together  yield
that
$$|\nabla |\psi|\,|\leq |(\nabla-iA)\psi|\leq\frac{C}\varepsilon,\quad
{\rm in}~\Omega\,.$$
Now, since
$$|\nabla u_\varepsilon|\leq \frac{C}\varepsilon\,$$
we deduce that
$$|\nabla|\varphi|\,|\leq \frac C{\varepsilon}\quad {\rm in}~\overline S\,.$$
Thus, the result of the lemma becomes a consequence of
\cite[Theorem~III.3]{BBH}.\hfill$\Box$\\

\section{Analysis of the Meissner state}\label{V-Sec-meissnerstate}
Let us recall the definition of $u_\varepsilon$ and
$C_0(\varepsilon)$ in Theorem~\ref{V-thm-kach3} and (\ref{V-C0})
respectively. This section is devoted to the analysis of the
following variational problem~:
\begin{equation}\label{V-EM}
M_0(\varepsilon,H)=\min_{A\in H^1(\Omega;\mathbb R^2)} \mathcal
G_{\varepsilon,H}(u_\varepsilon,A)\,.
\end{equation}
Since the function $u_\varepsilon$ is real-valued, one gets, for any
vector field $A$, the following decomposition~:
$$\mathcal G_{\varepsilon,H}(u_\varepsilon,A)=
C_0(\varepsilon)+\int_\Omega \left(|Au_\varepsilon|^2+ |{\rm
curl}\,A-H|^2\right)\,\md x.$$ Putting further
$$A=H\,\mathcal A,$$
\begin{equation}\label{V-J0}
J_0(\varepsilon)=\inf_{\mathcal A\in H^1(\Omega;\mathbb R^2)}\left[
\int_\Omega\left(|\mathcal A\,u_\varepsilon|^2+|{\rm curl}\,\mathcal
A-1|^2 \right)\,\md x\right],\end{equation} we get that
$$M_0(\varepsilon,H)=\inf_{A\in H^1(\Omega;\mathbb R^2)}
\mathcal G_{\varepsilon,H}(u_\varepsilon,A)=C_0(\varepsilon)+H^2J_0(\varepsilon),$$
and we are reduced to the analysis of the variational problem
(\ref{V-J0}).\\
Starting from a minimizing sequence (cf. \cite{SaSe}), it is
standard to prove that a minimizer $A_\varepsilon$ of (\ref{V-J0})
exists and satisfies the Coulomb gauge condition:
$${\rm div}\,A_\varepsilon=0\quad{\rm in}~\Omega,\quad n\cdot
 A_\varepsilon=0\quad{\rm on}~\partial\Omega,$$
where $n$ is the unit outward normal vector of the boundary of $\Omega$.\\
Notice also that $A_\varepsilon$ satisfies the Euler-Lagrange
equations~:
\begin{equation}\label{V-J0-EulerEq}
\nabla^\bot{\rm
curl\,}A_\varepsilon=u_\varepsilon^2\,A_\varepsilon\quad {\rm
in}~\Omega,\quad {\rm curl}\,A_\varepsilon=1\quad{\rm
on}~\partial\Omega.
\end{equation}
Here $\nabla^\bot=(-\partial_{x_2},\partial_{x_1})$
is the {\it Hodge gradient}.\\
Putting $h_\varepsilon={\rm curl}\,A_\varepsilon$, we get from the
first equation in (\ref{V-J0-EulerEq}) that
$A_\varepsilon=\frac1{u_\varepsilon^2}\nabla^\bot h_\varepsilon$. We
get also that $h_\varepsilon$ satisfies the equation:
\begin{equation}\label{V-hepsilon'}
-{\rm div}\left(\frac1{u_\varepsilon^2} \nabla
h_\varepsilon\right)+h_\varepsilon=0\quad{\rm in}~\Omega,\quad
h_\varepsilon=1\quad {\rm on}~\partial\Omega.
\end{equation}

\begin{lem}\label{V-lem-hepsilon}
The function $h_\varepsilon$ satisfies $0<h_\varepsilon<1$ in
$\Omega$, and it is the only function solving
(\ref{V-hepsilon'}).\\
Moreover, there exist constants $c_0,\varepsilon_0\in]0,1[$ such
that,
\begin{equation}\label{V-Eq-normhepsilon}
c_0\leq\|h_\varepsilon-1\|_{L^\infty(\Omega)}< 1,\quad
\forall~\varepsilon\in]0,\varepsilon_0].
\end{equation}
\end{lem}
\paragraph{\bf Proof.}
The property that $0<h_\varepsilon<1$ and the uniqueness of
$h_\varepsilon$ are direct applications of the Strong Maximum
Principle.\\
Let us now prove (\ref{V-Eq-normhepsilon}).  Assume by contradiction that there exists a
sequence converging to $0$, still denoted by $\varepsilon$, such
that
\begin{equation}\label{V-contradiction}
\lim_{\varepsilon\to0}\|h_\varepsilon-1\|_{L^\infty(\Omega)}=0.
\end{equation}
Let us take a compact subset $K\subset S$ (independent of
$\varepsilon$). Due to the asymptotic behaviour of $u_\varepsilon$
(it remains exponentially close to $1$ in $K$, see
Theorem~\ref{V-thm-kach3} and \cite[Proposition~5.1]{kach3}), it
results from (\ref{V-hepsilon'}) that $h_\varepsilon$ is bounded in
the $C^2$-norm of $K$. Thus, one can extract a subsequence of
$h_\varepsilon$, still denoted by $h_\varepsilon$, that converges to
a function $h\in C^2(K)$. The function $h$ satisfies the limiting
equation,
$$-\Delta h+h=0\quad{\rm in}~K.$$
Coming back to (\ref{V-contradiction}), we get that $h\equiv1$ in
$K$, hence not a solution of the limiting equation. Therefore, the
assertion in (\ref{V-Eq-normhepsilon}) holds.\hfill$\Box$\\

The next results
concern the case of our particular domains in (\ref{V-N}) and
(\ref{V-S}).

\begin{lem}\label{V-lem-disc}
With the assumptions (\ref{V-N}) and (\ref{V-S}), the function
$h_\varepsilon$ is radial, i.e. $h_\varepsilon(x)=\widetilde
h_\varepsilon(|x|)$, with $\widetilde h_\varepsilon$ being an
increasing function.
\end{lem}
\paragraph{\bf Proof.}
That $h_\varepsilon$ is radial follows by the uniqueness of the
solution of (\ref{V-hepsilon'}) and by the fact that $u_\varepsilon$
is also radial.\\
The solution $h_\varepsilon$ being radial, i.e.
$$h_\varepsilon(x)=\widetilde h_\varepsilon(|x|),\quad
\forall~x\in\Omega,$$ let us  show  that the function $\widetilde
h_\varepsilon$ is increasing. For simplicity of notation, we shall
remove the tilde and write $h_\varepsilon$ for $\widetilde
h_\varepsilon$. Notice that $h_\varepsilon$ satisfies the
differential equation~:
\begin{equation}\label{V-Eq-radialh}
\left\{\begin{array}{l} - h''_\varepsilon(r)-\displaystyle\frac1{r}
h'_\varepsilon(r)
+2\frac{u'_\varepsilon(r)}{u_\varepsilon(r)}\,h_\varepsilon'(r)+u_\varepsilon^2(r)\,h_\varepsilon(r)=0,\quad
0<r<1,\\
h_\varepsilon'(0)=0,\quad h_\varepsilon(1)=1.
\end{array}\right.
\end{equation}
Let us calculate $h''_\varepsilon(0)$. Since $h'_\varepsilon(0)=0$,
we have  $h''_\varepsilon(0)=\displaystyle\lim_{r\to0}
\displaystyle\frac{h'_\varepsilon(r)}{r}$. Substituting in
(\ref{V-Eq-radialh}), we get that
\begin{equation}\label{V-h''(0)}
h_\varepsilon''(0)=\frac12u_\varepsilon^2(0)\,h_\varepsilon(0)>0.
\end{equation}
Let us introduce the even extension of $h_\varepsilon$, namely the
function
$$f_\varepsilon(r)=\left\{\begin{array}{l}
h_\varepsilon(r)\quad(r>0),\\
h_\varepsilon(-r)\quad (r<0).
\end{array}\right.$$
Then $f_\varepsilon$ satisfies the equation,
\begin{equation}\label{V-Eq-fepsilon}
-f''_\varepsilon(r)-\displaystyle\frac1{|r|} f'_\varepsilon(r)
+2\frac{\widetilde u'_\varepsilon(r)}{\widetilde
  u_\varepsilon(r)}\,f_\varepsilon'(r)+\widetilde
u_\varepsilon^2(r)\,f_\varepsilon(r)=0,\quad
r\in]-r_2,r_2[\setminus\{0\},
\end{equation}
and it attains a local minimum at $0$. We emphasize also here that
$\widetilde u_\varepsilon$ denotes the even extension of
$u_\varepsilon$.\\
If $r_0\in]-1,1[$ (with $r_0\not=0$) is a critical point of
$f_\varepsilon$, then it follows from (\ref{V-Eq-fepsilon}) that~:
$$f''_\varepsilon(r_0)=\widetilde u_\varepsilon^2(r_0)\,f_\varepsilon(r_0)>0.$$
If $r_0=0$, the conclusion $f''_\varepsilon(0)>0$ still holds,
thanks
to (\ref{V-h''(0)}).\\
Now these observations lead to the conclusion that $f_\varepsilon$
attains its minimum at a unique point, and that this point is the
only  critical point for $f_\varepsilon$. As we know that
$f'_\varepsilon(0)=0$, we get that $f_\varepsilon$ attains its
minimum at $0$ and that it is increasing
in $[0,1[$. This achieves the proof of the lemma.\hfill$\Box$\\

The next lemma plays a distinguished role in the control of the
minimizing energy of `vortex balls'.

\begin{lem}\label{V-lem-controlgrad}
There exist constants $C>0$ and $\varepsilon_0>0$ such that
\begin{equation}\label{V-Eq-controlgrad}
\left\|\frac1{u_\varepsilon^2}\nabla
h_\varepsilon\right\|_{L^\infty(\Omega)}\leq 1, \quad
\forall~\varepsilon\in]0,\varepsilon_0].
\end{equation}
\end{lem}

\subsubsection*{ Proof.} Notice that by Lemma~\ref{V-lem-disc},
$h_\varepsilon$ is radial. Then the equation for $h_\varepsilon$ can
be written in the form:
$$-\left(\frac{h'_\varepsilon}{u_\varepsilon^2}\right)'(r)
-\frac1r\,\frac{h_\varepsilon'}{u_\varepsilon^2}(r)
+h_\varepsilon(r)=0,\quad\forall~r\in]0,1[.
$$
Integrating this equation between $0$ and $r\in]0,1[$ and using the
fact that $h_\varepsilon$ is increasing, $h_\varepsilon'\geq0$, we
obtain:
$$\left(\frac{h'_\varepsilon}{u_\varepsilon^2}\right)(r)\leq
\int_{0}^rh_\varepsilon(\widetilde r)\,\md \widetilde r\leq
r\|h_\varepsilon\|_{L^\infty([0,1])}\leq 1,$$ which is the result of
the lemma.
\hfill$\Box$\\

\section{Lower bound of the energy}\label{Section-Lowerbound}

\subsection{Construction of vortex-balls}\label{V-SubSec-VortexBalls}

We borrow some notation used in \cite{SaSe}. For a set $U\subset
\mathbb R^2$ we denote by $r(U)$ the {\it radius} of $U$, that is
the infimum over all finite coverings of $U$ by open balls
$B_1,B_2,\dots, B_k$ of the sum $r_1+r_2+\dots+r_k$. The important
property is that:
$$r(U)\leq \frac12\mathcal H^1(\partial U),$$ where $\mathcal H^1$ is
the one-dimensional Hausdorff measure.\\

\noindent From now on, we shall always work under the following
hypothesis:
$${\rm (H)}\quad
\exists\,c>0,~\forall~\varepsilon\in]0,1],\quad
\ell(\varepsilon)\leq c\ln|\ln\varepsilon|,\quad H\leq
c\,m_\varepsilon\,|\ln\varepsilon|\,,$$ where
\begin{equation}\label{m-epsilon}
m_\varepsilon=\inf_{x\in\overline\Omega}u_\varepsilon^2(x)\,.
\end{equation}
Notice that when $\ell(\varepsilon)=d\,\varepsilon$, $m_\varepsilon$
and $m_\varepsilon^{-1}$ are bounded. When $\ell(\varepsilon)\gg
\varepsilon$, we have by Theorem~\ref{L-infty-estimate},
$$m_\varepsilon=
A(a)\exp\left(-\frac{2\sqrt{a}\,\ell(\varepsilon)}\varepsilon\right)(1+o(1)),\quad
(\varepsilon\to0),
$$
where $A(a)>0$ is an explicit constant.

\begin{lem}\label{V-lemSaSeN}
Let $(\psi,A)$ be a minimizer of (\ref{V-EGL}) and
$\varphi=\displaystyle\frac{\psi}{u_\varepsilon}$. Then, under the
hypotheses {\rm (H)} above,
 there exist constants $C>0$ and $\varepsilon_0>0$ such that,
for all
$$\delta\in
\left]\left(|\ln\varepsilon|+|\ln m_\varepsilon|\right)
\sqrt{\frac{\sqrt{\varepsilon}}{\varepsilon_0}}\,,1\right[\quad {\rm
and}~\varepsilon\in]0,\varepsilon_0],$$ we have
\begin{equation}\label{V-Eq-SaSe1-N}
r\left(\{x\in \overline \Omega~:~|\varphi(x)|\leq
1-\delta\}\left)\,\leq C\,\frac{\sqrt{\varepsilon}
|\ln\varepsilon|^2}{\delta^2}\right.\right..
\end{equation}
\end{lem}
\paragraph{\bf Proof.}
We have the following decomposition of the energy,
$$\mathcal G_{\varepsilon,H}(\psi,A)=C_0(\varepsilon)
+\mathcal F_{\varepsilon,H}(\varphi,A),
$$
where the functional $\mathcal F_{\varepsilon,H}(\varphi,A)$ has
been introduced in (\ref{V-reducedfunctional}).\\
Using $\left(u_\varepsilon,\frac1{u_\varepsilon^2}\nabla^\bot
h_\varepsilon\right)$ as a test configuration for the functional
(\ref{V-EGL}), we get
\begin{equation}\label{M-UB}
\mathcal F_{\varepsilon,H}(\psi,A)\leq \widetilde
c\,H^2.\end{equation} We infer from (\ref{M-UB}),
$$\int_\Omega|(\nabla-iA)\varphi|^2\,\md x
+\frac{1}{2\varepsilon^2}\int_S(1-|\varphi|^2)^2\,\md
x\leq2\,\widehat c\,m_\varepsilon^{-2}H^2\,,$$ where $m_\varepsilon$
is introduced in (\ref{m-epsilon}).\\
By Lemma~\ref{V-lem-psi<u},
$|\varphi|\leq 1$, hence
\begin{equation}\label{Pinning-N} \int_N(1-|\varphi|^2)^2\,\md x\leq
2|N|\leq 4\pi\ell(\varepsilon)\leq
C\,\varepsilon|\ln m_\varepsilon|\,.\end{equation} Now, by
(\ref{Pinning-N}) and the diamagnetic inequality,
$|(\nabla-iA)\varphi|\geq |\nabla |\varphi|\,|$,  we deduce
\begin{equation}\label{SaSe-upperbound}
\int_\Omega|\nabla|\varphi|\,|^2\,\md x
+\frac{1}{2\widetilde{\varepsilon}^2}\int_\Omega(1-|\varphi|^2)^2\,\md
x\leq M\,,\end{equation} where
$$\widetilde\varepsilon=\sqrt{\varepsilon}\,,\quad M=2\,\widetilde
c\left(m_\varepsilon^{-2}H^2+|\ln m_\varepsilon|\right).$$ Since
$\varphi$ is a $H^2$-function, then   it can be approximated in
 the
$L^\infty$-norm  by means of smooth functions.  This permits us to
conclude (\ref{V-Eq-SaSe1-N})
from Proposition~\ref{Sandier1}.\hfill$\Box$\\

The next proposition provides us, as in \cite{SaSe}, with the
construction of suitable `vortex-balls'.

\begin{prop}\label{V-lem-vortexballs}
Let $(\psi,A)$ be a minimizer of (\ref{V-EGL}) and
$\varphi=\displaystyle\frac{\psi}{u_\varepsilon}$. Then, under the
hypotheses {\rm (H)}, for each $p\in]1,2[$, $\alpha\in]0,1/2[$,
$n\in\mathbb N$, 
there exist constants $C>0$,
$\gamma\in]0,1/2[$,  and for each  $\varepsilon^\alpha\ll\eta\ll1$,
there exists a family of disjoint balls
$\{B((a_i,r_i)\}_i$ satisfying the following properties:
\begin{enumerate}
\item $w=\{x\in \overline\Omega~:~|\varphi(x)|\leq
  1-|\ln\varepsilon|^{-n}\}
\subset\displaystyle\cup_{i}B(a_i,r_i)$.
\item $\displaystyle\sum_{i}r_i\leq
\eta$.
\item Letting $d_i$ be the degree of the
function $\varphi/|\varphi|$ restricted to $\partial B(a_i,r_i)$ if
$B(a_i,r_i)\subset\Omega$ and $d_i=0$
  otherwise, then we have:
\begin{eqnarray}\label{V-Eq-LBest}
&&\hskip-0.5cm \int_{B(a_i,r_i)\setminus\omega}
u_\varepsilon^2|(\nabla-iA)\varphi|^2\,\md x+\int_{B(a_i,r_i)}|{\rm
  curl}\,A-H|^2\,\md x\geq\\
&&\hskip3.5cm 2 \pi|d_i|
\left(\min_{B(a_i,r_i)}u_\varepsilon^2\right)
\left(\ln\frac{\eta}{\varepsilon^\alpha}-C|\ln\varepsilon|^{-n}\right).
\nonumber
\end{eqnarray}
\item
$\left\|2\pi\displaystyle\sum_i d_i\delta_{a_i}-{\rm
    curl}\big{(}A+(i\varphi,\nabla_A\varphi)\big{)}
\right\|_{W^{-1,p}_0(\Omega)}\leq
\max(|\ln\varepsilon|^{2-2n},\eta^\gamma).$
\end{enumerate}
\end{prop}
\paragraph{\bf Proof.}
Let us take  $\delta=|\ln\varepsilon|^{-n}$ in the statement of
Lemma~\ref{V-lemSaSeN}. We emphasize that under the hypothesis (H),
our choice of $\delta$ verifies the hypothesis of
Lemma~\ref{V-lemSaSeN}, namely, $\delta\in
\left]\left(|\ln\varepsilon|+|\ln m_\varepsilon|\right)
\sqrt{\frac{\sqrt{\varepsilon}}{\varepsilon_0}}\,,1\right[$.\\
Points (1) and (2) of Proposition~\ref{V-lem-vortexballs}  are now
direct consequences of Lemma~\ref{V-lemSaSeN}
and the first point of Proposition~\ref{Sandier2}.\\
Let us prove now Point (3). By the estimate on $r(w)$, we get
for a given   $\alpha\in]0,\frac12[$, $r(w)<\varepsilon^\alpha$
provided that $\varepsilon$ is
sufficiently small.\\
Now, notice that in $B(a_i,r_i)\setminus w$,
$$
|(\nabla-iA)\varphi|^2\geq |\varphi|^2\left|(\nabla-iA)
\frac{\varphi}{|\varphi|}\right|^2,\quad
\frac{|\varphi|^2-1}{|\varphi|^2}\geq -2|\ln\varepsilon|^{-n}.$$
Hence, we get the desired conclusion by applying the third point of
Proposition~\ref{Sandier2} to the function $\varphi/|\varphi|$.\\
Notice that the function $\varphi$ and the balls $B(a_i,r_i)$
satisfy the hypotheses of Proposition~\ref{Sandier3} with
$M=\mathcal O( |\ln\varepsilon|^2)$. The application of this
proposition yields the conclusion in the last point of
Proposition~\ref{V-lem-vortexballs}.
\hfill$\Box$\\

\begin{rem}
When performing the previous argument with $\delta=\varepsilon^\beta$
and $\beta\in]0,\frac12[$ sufficiently small, we get improved remainders
  in (3)-(4) of Proposition~\ref{V-lem-vortexballs}, a power of $\varepsilon$,
  but valid for larger values of $\eta$. This permits to treat the
  case when the thickness 
$\ell(\varepsilon)\leq c\,\varepsilon|\ln\varepsilon|$ and $c\in]0,1[$
  is sufficiently small.
\end{rem}

We follow the usual terminology and call the balls constructed in
Proposition~\ref{V-lem-vortexballs} `vortex-balls'.\\

We conclude with the following theorem.

\begin{thm}\label{LowerBound-thm}
Let $(\psi,A)$ be a minimizer of (\ref{V-EGL}) and
$\varphi=\displaystyle\frac{\psi}{u_\varepsilon}$. Then, under the
hypothesis {\rm (H)}, for each
 $\alpha\in]0,1/2[$ and $n\in\mathbb N$,
there exist a constant $C>0$ and a family of disjoint balls
$\{B((a_i,r_i)\}_i$ such that~:
\begin{enumerate}
\item $\displaystyle\sum_i r_i\leq
  C\varepsilon^{\alpha'}$, $(\alpha'\in]0,\frac12-\alpha[)$;\\
\item $|\varphi|\geq\frac12$ on $\Omega\setminus \cup_iB(a_i,r_i)$.
\item Letting $d_i$ be the degree of the
function $\varphi/|\varphi|$ restricted to $\partial B(a_i,r_i)$ if
$B(a_i,r_i)\subset\Omega$ and $d_i=0$
  otherwise, then  we have:
\begin{equation}\label{V-Eq-LBest'}
\begin{split}
\mathcal F_{\varepsilon,H}(\varphi,A)\,\geq&
\,H^2J_0(\varepsilon)\\
&+2\pi \sum_{d_i\geq0}
\left[\alpha\left(\min_{B(a_i,r_i)}u_\varepsilon^2\right)|\ln\varepsilon|-
2 H\right]d_i-CH|\ln\varepsilon|^{-n}\,;
\end{split}
\end{equation}
\item $\sum_i|d_i|\leq C\,m_\varepsilon |\ln\varepsilon|$\,.
\end{enumerate}
\end{thm}
\paragraph{\bf Proof}
Applying Proposition~\ref{V-lem-vortexballs}  with
$\eta=\varepsilon^{\alpha'}$ and with 
$\alpha$ replaced by $\beta\in]0,\frac12[$ to be chosen sufficiently
  close to $\frac12$, we get a family of balls satisfying in
  particular the first two  assertions of  the theorem. 
\subsubsection*{\it Total degree}\ \\
We start by proving an upper bound on the total degree
$\sum_i|d_i|$. The starting point is by noticing that there exists a
constant $c>0$ such that
$$\mathcal F_{\varepsilon,H}(\varphi,A)\leq c\,H^2.$$
Then, by applying Point (3) of Proposition~\ref{V-lem-vortexballs}
with $\eta=\varepsilon^{\alpha'}$, we deduce the existence of a constant $\widetilde
c>0$ such that
$$2\pi\,m_\varepsilon|\ln\varepsilon|\sum_i|d_i|\leq \widetilde
c\,H^2,$$ hence Point (4) of the theorem is proved.
\subsubsection*{\it A rough lower bound of the energy}\ \\
We put
$$\widetilde \Omega=\Omega\,\Big{\backslash}\bigcup_
{B(a_i,r_i)\subset \Omega}B(a_i,r_i),$$
$$j=\left(i\varphi,\nabla_{A}\varphi\right),\quad j'=j-\frac{H}{u_\varepsilon^2}\nabla^\bot h_\varepsilon,\quad
A'=A-\frac{H}{u_\varepsilon^2}\nabla^\bot h_\varepsilon.$$ Since
$|\varphi|\leq 1$, then $$|j|\leq |\varphi|\,|\nabla_A\varphi|\leq
|\nabla_A\varphi|,$$ and consequently,  we have
\begin{align*}
\begin{split}
\mathcal
F_{\varepsilon,H}(\varphi,A,\widetilde\Omega)\geq&\int_{\widetilde\Omega}\left(
u_\varepsilon^2|j|^2+|{\rm curl}\,A-H|^2\right)\md x .
\end{split}
\end{align*}
Now a direct calculation yields,
\begin{align*}
\begin{split}
\mathcal F_{\varepsilon,H}(\varphi,A,\widetilde\Omega)\geq&
H^2J_0(\varepsilon,\widetilde\Omega)+ 2H\int_{\Omega} (h_\varepsilon-1)\big{[}{\rm curl}( A'+j')\big{]}\,\md x\\
&-2H\int_{\cup_iB(a_i,r_i)}\left[(h_\varepsilon-1)\,{\rm
curl}\,A'-j'\cdot\nabla^\bot h_\varepsilon\right]\,\md x,
\end{split}
\end{align*}
where
$$J_0(\varepsilon,\widetilde\Omega)=\int_{\widetilde\Omega}
\left(\frac1{u_\varepsilon^2}|\nabla
u_\varepsilon|^2+|h_\varepsilon-1|^2\right)\,\md x.$$ Using
Lemma~\ref{V-lem-controlgrad} and the Cauchy-Schwarz inequality, it
is easy to prove that
$$
\left|\int_{\cup_iB(a_i,r_i)}\left[(h_\varepsilon-1)\,{\rm
curl}\,A'-j'\cdot\nabla^\bot h_\varepsilon\right]\,\md x\right|\leq
C\,H\sum_ir_i,$$ and
$$\left|J_0(\varepsilon)-J_0(\varepsilon,\widetilde\Omega)\right|\leq C\sum_{i}r_i.$$
Therefore, we obtain
\begin{equation}\label{V-lowerbound-exterior}
\mathcal F_{\varepsilon,H}(\varphi,A,\widetilde\Omega)\geq
H^2J_0(\varepsilon)+ 2H\int_{\Omega} (h_\varepsilon-1)\big{[}{\rm
curl}( A'+j')\big{]}\,\md x-CH^2\sum_ir_i.
\end{equation}
Since ${\rm curl}( A'+j')={\rm curl}\left(
A+(i\varphi,\nabla_A\varphi)\right)$, then by Point (4) in
Proposition~\ref{V-lem-vortexballs}, we rewrite the above lower
bound in the form
\begin{equation}\label{V-lowerbound-exterior}
\mathcal F_{\varepsilon,H}(\varphi,A,\widetilde\Omega)\geq
H^2J_0(\varepsilon)+ 4\pi H\sum_i
d_i(h_\varepsilon-1)(a_i)-CH|\ln\varepsilon|^{-n}.
\end{equation}
We have also by Point (3) of Proposition~\ref{V-lem-vortexballs},
$\forall \alpha\in]0,\frac12[$ and when $\varepsilon$ is sufficiently small,
\begin{equation}\label{V-lowerbound-interior}
\sum_i\mathcal
F_{\varepsilon,H}\left(\varphi,A,B(a_i,r_i)\right)\geq 2\pi\alpha \sum_i
|d_i|\left(\min_{B(a_i,r_i)}u_\varepsilon^2\right)|\ln\varepsilon|.
\end{equation}
Therefore, we obtain the lower bound
\begin{eqnarray}\label{V-lowerbound-global}
\mathcal F_{\varepsilon,H}(\varphi,A)&\geq&
H^2J_0(\varepsilon)+2\pi\alpha\,|\ln\varepsilon| \sum_i
|d_i|\left(\min_{B(a_i,r_i)}u_\varepsilon^2\right)\\&&+ 4\pi H\sum_i
d_i(h_\varepsilon-1)(a_i)-CH|\ln\varepsilon|^{-n}.\nonumber
\end{eqnarray}
Since $0<h_\varepsilon<1$, the above lower bound is sufficient to
deduce (\ref{V-Eq-LBest'}).
\hfill$\Box$\\

\section{Proofs of main results}\label{proofs}
\subsection{A vortex-less regime}\label{existence-lambda}

Let us recall the definition of the constant $m_\varepsilon$
introduced in (\ref{m-epsilon}). We recall also that $(\psi,A)$
always denotes a minimizer of (\ref{V-EGL}) and that
$\varphi=\frac{\psi}{u_\varepsilon}$.\\

The aim of this subsection is to prove the following
theorem.  

\begin{thm}\label{existence-lambda}
There exists a constant $\lambda>0$ such that if $\ell\leq
\mathcal O(\varepsilon\ln|\ln\varepsilon|)$ and if
the magnetic field satisfies
$$H\leq \lambda m_\varepsilon |\ln\varepsilon|\,,$$
then
$$|\varphi|\geq \frac12\quad {\rm in}~ \overline S\,,$$
and we have the energy estimate as $\varepsilon\to0$
$$\int_\Omega\left(|(\nabla-iA')\varphi|^2+
|{\rm curl}\,A-H\,h_\varepsilon|^2\right)\,\md x+\frac1{\varepsilon^2}
\int_S(1-|\varphi|^2)^2\,\md x\ll m_\varepsilon^4\,.$$
Here
$$A'=A-\frac{H}{u_\varepsilon^2}\nabla^\bot h_\varepsilon$$
and $h_\varepsilon$ is the function introduced in (\ref{V-hepsilon'}).
\end{thm}

One essential step towards the proof of Theorem~\ref{existence-lambda}
is a further useful splitting of the
energy due to B\'ethuel-Rivi\`ere (cf. \cite{BeRi}).

\begin{lem}\label{Beth-Riv}
Consider $(u,A)\in H^1(\Omega;\mathbb C)\times H^1(\Omega;\mathbb
R^2)$ and define
$$A'=A-\frac{H}{u_\varepsilon^2}\nabla^\bot h_\varepsilon,$$
where $u_\varepsilon$ and $h_\varepsilon$ are introduced in
Theorem~\ref{V-thm-kach3} and (\ref{V-hepsilon'}) respectively. Then
we have the decomposition of the energy,
\begin{align*}
\begin{split}
\mathcal F_{\varepsilon,H}(u,A)=&H^2J_0(\varepsilon)
+\int_\Omega\left( u_\varepsilon^2|(\nabla-iA')u|^2+|{\rm
curl}\,A'|^2\right)\md x
\\
&+\frac1{\varepsilon^2}\int_S u_\varepsilon^4(1-|u|^2)^2\md x+
2H\int_\Omega (h_\varepsilon-1)\bigg{[}{\rm curl}\big{(}
A'+(iu,\nabla_{A'}u)\big{)}\bigg{]}\md x\\
&+H^2\int_\Omega\frac1{u_\varepsilon^2}\left(|u|^2-1\right)
|\nabla h_\varepsilon|^2\,\md x.
\end{split}
\end{align*}
Here, the functional $\mathcal F_{\varepsilon,H}$ and the energy
$J_0(\varepsilon)$ are introduced in (\ref{V-reducedfunctional}) and
(\ref{V-J0}) respectively. \end{lem}
\break
\subsubsection*{Proof of Theorem~\ref{existence-lambda}}
\paragraph{\it Existence of $\lambda$\,.}\ \\
Let us choose $\lambda>0$ in such a way that when the magnetic field
satisfies
$$H\leq \lambda m_\varepsilon|\ln\varepsilon|\,,
$$
all
the degrees $d_i$ given by Theorem~\ref{LowerBound-thm} are equal
to zero, so that the energy of a minimizer becomes close to that of
the Meissner state.\\
Matching the lower bound (\ref{V-Eq-LBest'}) with the upper
bound
\begin{equation}\label{Upperbound-trivial}
\mathcal F_{\varepsilon,H}(\varphi,A)\leq
H^2J_0(\varepsilon)\,,
\end{equation} we obtain,
$$\sum_{d_i>0}
\left[\alpha\,m_\varepsilon|\ln\varepsilon|-2
H\right]d_i+CH\sum_{d_i<0}|d_i|-CH|\ln\varepsilon|^{-n}\leq 0\,.$$
By our choice of
$H$ and $\ell(\varepsilon)$, we get that for any $n\in\mathbb N$,
there exists a constant $C>0$ such that
$$\sum_{d_i>0}
\left[\alpha\, m_\varepsilon|\ln\varepsilon|-2
H\right]d_i+CH\sum_{d_i<0}|d_i|\leq |\ln\varepsilon|^{-n}.$$ Now
we get by choosing $\lambda<2\alpha$ that, for all
$\lambda\leq\lambda_1$,
$$0\leq H\min\left(\frac{\alpha}{\lambda}-2,C\right)\sum_i|d_i|\leq
 |\ln\varepsilon|^{-n}.$$
Therefore, for all $i$, $d_i=0$ and there does not exist vortices.\\
Moreover, the lower bound (\ref{V-Eq-LBest'}) becomes when
$\lambda\leq \lambda_1$,
\begin{equation}\label{Energy-vortexless}
\mathcal F_{\varepsilon,H}(\varphi,A)\geq
H^2J_0(\varepsilon)-C|\ln\varepsilon|^{-n},\end{equation}
hence when combined
with the upper bound (\ref{Upperbound-trivial}) together with the
energy expansion of Lemma~\ref{Beth-Riv} and Item (4) of
Proposition~\ref{V-lem-vortexballs}, we are able to
deduce the energy estimate of Theorem~\ref{existence-lambda}.\\
\paragraph{\it Proof of $|\varphi|\geq\frac12$.}\ \\
Now, as we have the energy estimate, we deduce that
$$\frac1{\varepsilon^2}
\int_S(1-|\varphi|^2)\,\md x\ll 1 \quad (\varepsilon\to0),$$
where we have also used that, in $\overline S$, $u_\varepsilon\geq c_1$
for some explicit constant $c_1>0$. Therefore, we infer from
Lemma~\ref{BBH-thm3.3} that $|\varphi|\geq \frac12$ in $\overline S$
\,.\hfill$\Box$\\

It results from Theorem~\ref{existence-lambda} a uniform estimate of $|A'|$.

\begin{corol}\label{Linfty-A'}
With the notations and hypothesis of Theorem~\ref{existence-lambda}, and
if the magnetic field satisfies
$$H\leq \lambda m_\varepsilon|\ln\varepsilon|\,,$$
then we have as $\varepsilon\to0$,
$$\|A'\|_{H^2(\Omega)}\ll m_\varepsilon^2\,.$$
\end{corol}
\paragraph{\bf Proof.}
Notice that since $u_\varepsilon$ and $h_\varepsilon$ are radial, then
$$\nabla u_\varepsilon (x)=n(x) u'_\varepsilon(|x|),\quad
\nabla h_\varepsilon (x)=n(x) h'_\varepsilon(|x|),$$
where $n(x)=\frac{x}{|x|}$ for all $x\in\mathbb R^2\setminus\{0\}$.\\
Therefore, the vector field $A'$ satisfies the properties inferred from $A$,
$${\rm div}\,A'=0\quad{\rm in}~\Omega,\quad
n(x)\cdot A'=0\quad{\rm on}~\partial\Omega\,.$$
Thus, it is a result of the Poincar\'e Lemma that
$$\|A'\|_{H^2(\Omega)}\leq C\|{\rm curl}\,A'\|_{L^2(\Omega)}\,$$
for some explicit geometric constant $C>0$.\\
Now, from the energy estimate of Theorem~\ref{existence-lambda}, we
conclude the result of the corollary.\hfill$\Box$\\

The next lemma is now an essential step in proving
Theorems~\ref{mainthm}-\ref{mainthm2}.

\begin{lem}\label{L2-Bndestimates}
With the notations and hypothesis of Theorem~\ref{existence-lambda}, and
if the magnetic field satisfies
$$H\leq \lambda m_\varepsilon|\ln\varepsilon|\,,$$
then we have as $\varepsilon\to0$,
$$\varepsilon\big{\|}\,n(x)\cdot(\nabla-iA')\varphi\,\big{\|}
_{L^2(\partial N)}\ll 1\,,$$
and
$$\big{\|}\,[\,\varphi\,]_N\big{\|}_{L^2(\mathbb S_R^1)}\ll 1\,.$$
Here  $n(x)=\displaystyle\frac{x}{|x|}$ for all $x\in\mathbb
R^2\setminus\{ 0\}$, and
$[\cdot]_N$ is the jump across $N$ introduced in (\ref{jump}).
\end{lem}
\paragraph{\bf Proof.}
Let us introduce the two domains
$$\mathcal S_{1}=D(0,R-\ell)\setminus D(0,R/2-\ell),\quad
\mathcal S_{2}=D(0,[R+1]/2+\ell)\setminus D(0,R+\ell)\,.$$
Notice that the domains $\mathcal S_1$ and $\mathcal S_2$ can be identified by those
corresponding to $\ell=0$ via the translations
$$T_\pm~: x=Re^{i\theta}\mapsto (R\pm\ell)e^{i\theta}\,.$$
In order to prove the first statement of the lemma,
we have only to establish (thanks to the trace
theorem),
\begin{equation}\label{CLAIM-100}
\sum_{i=1}^2
\varepsilon\left\|\,|n(x)\cdot(\nabla-iA')\varphi|\,\right\|_{H^1(\mathcal
  S_i)}
\ll1\,\quad(\varepsilon\to0).\end{equation}
Notice that, by the energy estimate of Theorem~\ref{existence-lambda},
we have only to estimate
$\|\nabla\,|n(x)\cdot(\nabla-iA')\varphi|\,\|_{L^2}$.\\
Notice that since $h_\varepsilon$ is radial,
$$n(x)\cdot\nabla^\bot h_\varepsilon=0,$$
hence
$$n(x)\cdot (\nabla-iA')\varphi=n(x)\cdot(\nabla-iA)\varphi\,.$$
Now, one deduces from Corollary~\ref{Linfty-A'} together with the
energy estimate of Theorem~\ref{existence-lambda},
\begin{equation}\label{CLAIM}
\sqrt{\varepsilon}\|A\|_{L^\infty(\Omega)}\ll m_\varepsilon,\quad
\sqrt{\varepsilon}\|A\|_{H^1(\Omega)}\ll m_\varepsilon,\quad
\sqrt{\varepsilon}\|\nabla\varphi\|_{L^2(\Omega)}\ll m_\varepsilon,
\end{equation}
and consequently, we obtain
$$\varepsilon\|A\cdot\nabla \varphi\|_{L^2(\Omega)}\ll m_\varepsilon,\quad
\varepsilon\|\,|A|^2\varphi\|_{L^2(\Omega)}\ll m_\varepsilon\,.$$
Now,
$$\Delta
\varphi=(\nabla-iA)^2\varphi+2iA\cdot\nabla\varphi+|A|^2\varphi\,$$
where
\begin{eqnarray*}
(\nabla-iA)^2\varphi&=&2\frac{\nabla u_\varepsilon}{u_\varepsilon}\cdot
(\nabla-iA)\varphi+(\nabla-iA)\cdot
  u_\varepsilon^2(\nabla-iA)\varphi\\
&=&2\frac{\nabla u_\varepsilon}{u_\varepsilon}\cdot
(\nabla-iA')\varphi+(\nabla-iA)\cdot
  u_\varepsilon^2(\nabla-iA)\varphi\\
&&
[\,{\rm since~}u_\varepsilon~{\rm is~ radial}\,]\,,\end{eqnarray*}
and
$$\varepsilon\left\|\frac{\nabla u_\varepsilon}{u_\varepsilon}
\right\|_{L^\infty(\Omega)}\leq m_\varepsilon^{-1}\,.$$
Therefore, we deduce from the energy estimate of
Theorem~\ref{existence-lambda} and the G-L equation of $\varphi$ that
$$\varepsilon\|\Delta\varphi\|_{L^2(\Omega)}\ll m_\varepsilon\,.$$
Now, by the standard  elliptic regularity theorem, we get
\begin{equation}\label{H2-estimate-phi}
\varepsilon\|\varphi\|_{H^2(\Omega)}\ll m_\varepsilon\,.\end{equation}
Therefore, we obtain,
$$\varepsilon\|(\nabla-iA)\varphi\|_{H^1(\Omega)}\ll 1\,.$$
As we have pointed out, this is now sufficient to deduce the first result of
the lemma.\\
We prove now the second statement. Let $(r,\theta)$ be polar
coordinates. Let us introduce the following function in $[0,2\pi[$,
$$f(\theta)=\varphi\big{(}(R+\ell)e^{i\theta}\big{)}
-\varphi\big{(}(R-\ell)e^{i\theta}\big{)}\,.$$
We claim that
\begin{equation}\label{f->0}
\forall~\theta\in[0,2\pi[,\quad
\lim_{\varepsilon\to0}|f(\theta)|=0\,.
\end{equation}
Actually, by the mean value theorem, we have
$$
|f(\theta)|\leq 2\ell\,
\|\partial_r\varphi(\cdot,\theta)\|_{L^\infty([R/2,(R+1)/2]}\,.
$$
On the other hand,
invoking (\ref{H2-estimate-phi}) and the sobolev embedding theorem, we
deduce that
$$|f(\theta)|\ll 2C_\theta
\frac{\ell(\varepsilon)}\varepsilon\,m_\varepsilon\,,$$
for a constant $C_\theta>0$.\\
If $\ell=\mathcal O(\varepsilon)$, this yields
(\ref{f->0}). Otherwise, if $\ell(\varepsilon)\gg \varepsilon$,
we have
$$\frac{\ell}\varepsilon m_\varepsilon\leq \frac{\ell}\varepsilon
\exp\left(\frac{-2\sqrt{a}\,\ell}\varepsilon\right)\ll 1,$$
which again yields  (\ref{f->0}).\\
Now, since $|\varphi|$ has values in $[0,1[$,
$\|f(R,\cdot)\|_{L^2([0,2\pi[)}$ is bounded. Therefore, we deduce by
the Lebesgue dominate convergence theorem that
$$\lim_{\varepsilon\to0}\|f(R,\cdot)\|_ {L^2([0,2\pi[)}=0,$$
which is nothing but the second statement of the lemma.\hfill$\Box$\\

\subsection{Proof of Theorems~\ref{mainthm}-\ref{mainthm2}}\ \\
Let us notice that
$$\left|n(x)\cdot\left[
\frac{(\nabla-iA)\psi}{\psi}-\frac{\nabla
  u_\varepsilon}{u_\varepsilon}\right]\,\right|
=\frac1{|\varphi|}\left|n(x)\cdot(\nabla-iA')\varphi\right|\,.$$
Thus, in the regime of Theorem~\ref{existence-lambda}, we have
$$
\varepsilon\left\|n(x)\cdot\left[
\frac{(\nabla-iA)\psi}{\psi}-\frac{\nabla
  u_\varepsilon}{u_\varepsilon}\right]\,\right\|_{L^2(\partial N)}
\leq 2\varepsilon
\left\|n(x)\cdot(\nabla-iA')\varphi\right\|_{L^2(\partial N)}\ll1\,,$$
where the last conclusion is due to Lemma~\ref{L2-Bndestimates}.\\
Invoking Theorem~\ref{mainthm-H=0}, we deduce the formulas
(\ref{jump-cond1}) and (\ref{jump-cond1'}).\\
Now, it remains to prove that
\begin{equation}\label{proof-psi-jump}
\big{\|}\,[\psi]_N\,\big{\|}_{L^2(\mathbb S_R^1)}\ll1\,.\end{equation}
But,
by Lemma~\ref{L2-Bndestimates}, it is easy to deduce
(\ref{proof-psi-jump}) since $\psi=u_\varepsilon\,\varphi$, and by
Lemma~\ref{mainthm-H=0}, $u_\varepsilon$ converges uniformly on
$\partial N$ to an
explicit constant $\mathcal A>0$.\\
Therefore, we get now that (\ref{jump-cond2}) and
(\ref{jump-cond2'}) hold.
This finally achieves the proof of Theorems~\ref{mainthm}~\&~\ref{mainthm2}.\hfill$\Box$\\

\subsection{Proof of Theorem~\ref{thm-current}}\ \\
Let $j=(i\psi,(\nabla-iA)\psi)$
and  $j_\varphi=(i\varphi,(\nabla-iA)\varphi)$, where
$\varphi=\frac{\psi}{u_\varepsilon}$.\\
Integrating by parts, we get
$$\int_{N}{\rm curl}\,j_\varphi=\int_{|x|=R+\ell}\tau(x)\cdot
j_\varphi-\int_{|x|=R-\ell}\tau(x)\cdot j_\varphi\,,$$
where $\tau(x)=\frac{x^\bot}{|x|}$ is the tangential vector of any
circle in $\mathbb R^2$.\\
Notice that $j_\varphi=\frac1{u_\varepsilon^2}\,j$. Then since the function
$u_\varepsilon$ is radial and by Theorem~\ref{mainthm-H=0},
$$\left\|u_\varepsilon- C\right\|_{L^\infty(\partial N)}\to0\quad({\rm as}~\varepsilon\to0)$$
for an explicit constant $C>0$, we deduce that
$$\left|
\int_{|x|=R+\ell}\tau(x)\cdot j
-\int_{|x|=R-\ell}\tau(x)\cdot j\right|= C\left|\int_N{\rm
  curl}\,j_\varphi\right|(1+o(1))\quad{\rm as}~\varepsilon\to0\,.$$
But, since we have no vortices, we get from Point (4) of
Proposition~\ref{V-lem-vortexballs}:
$$\int_N{\rm curl}\,j_\varphi=-\int_N{\rm curl}\,A+o(1)\quad
(\varepsilon\to0)\,.$$
Now by the energy estimate of Theorem~\ref{existence-lambda}, we
deduce
\begin{eqnarray*}
\left|\int_N{\rm curl}\,A\right|&\leq&|N|^{1/2}\,\|{\rm
  curl}\,A\|_{L^2(\Omega)}\\
&\lesssim&\varepsilon^{1/2}|\ln\varepsilon|.
\end{eqnarray*}
This achieves the proof of Theorem~\ref{thm-current}.\hfill$\Box$\\

\section*{Acknowledgements}
The author would like to thank B. Helffer for the interest he
owed to this work and for his many valuable suggestions. He would also
like to thank E. Sandier for  fruitful discussions. 
This work has been supported by the European Research Network
  `Post-doctoral Training Program in Mathematical Analysis of Large
  Quantum Systems' with contract number HPRN-CT-2002-00277, and the ESF
  Scientific Programme in Spectral Theory and Partial Differential
  Equations (SPECT).

\appendix

\section{The canonical equation}\label{Appendix-Canon}
Let $\beta=\beta(a,d)>e^{-\sqrt{2}\,d}$
and $A=A(a,d)>0$ be the solutions of the
following equations:
$$\left\{
\begin{array}{l}
\displaystyle\frac{2\sqrt{2}\,\beta\,e^{\sqrt{2}\,d}}{(\beta\,e^{\sqrt{2}\,d}+1)^2}
=
\sqrt{a}\left[e^{\sqrt{a}\,d}-
e^{-\sqrt{a}\,d}\right]A\,,\\
\displaystyle\frac{\beta\,e^{\sqrt{2}\,d}-1}{\beta\, e^{\sqrt{2}\,d}+1}
=\left[e^{\sqrt{a}\,d}+
e^{-\sqrt{a}\,d}\right]A\,.
\end{array}\right.$$
Putting $b=e^{\sqrt{a}\,d}+
e^{-\sqrt{a}\,d}$ and $c=\sqrt{a}\left[e^{\sqrt{a}\,d}-
e^{-\sqrt{a}\,d}\right]$, the solution $(\beta,A)$ can be expressed
explicitly by,
\begin{eqnarray}
\beta(a,d)&=&\sqrt{2}\left(\frac
bc+\sqrt{\left(\frac{b}{c}\right)^2+\frac12}\,\right)
e^{-\sqrt{2}\,d}\label{beta}\\
A(a,d)&=&
\frac{\beta(a,d)\,e^{\sqrt{2}\,d}-1}
{b\left(\beta(a,d)\,e^{\sqrt{2}\,d}+1\right)}\label{A(a,d)}\,.
\end{eqnarray}
With this choice, it is easily  checked
that the function $U$ given by (\ref{CanonEq1}) is
a solution of Eq. (\ref{CanonEq}). The aim of this appendix is to show
that (\ref{CanonEq1}) is the only bounded and positive solution of
(\ref{CanonEq}).

\begin{thm}\label{CanonEq-thm}
In the class of functions $\mathcal C=\{u\in H^1_{\rm loc}(\mathbb
R^2)\cap L^\infty(\mathbb R^2)~:~u\geq 0\}$, Eq. (\ref{CanonEq})
admits a unique non-trivial solution 
$$\mathbb R^2\ni(x_1,x_2)\mapsto U(x_2)$$
given by (\ref{CanonEq1}).
\end{thm}
\paragraph{\bf Proof.}
Since the proof is very close to that of \cite[Theorem~1.5]{kach3},
we sketch only the main steps.\\
By adjusting the proof of \cite[Lemma~4.2]{kach3}, we obtain that if
$u\not\equiv0$ solves (\ref{CanonEq}), then $0<u<1$ in $\mathbb
R^2$. This permits us, when following step by step the proof of
\cite[Lemma~4.3]{kach3} and \cite[Lemma~5.3]{LuPa96},
to get a positive constant $C\in]0,1[$
such that for any solution $u$ of (\ref{CanonEq}) in $\mathcal C$, we
have
\begin{equation}\label{kach3-4.9}
\inf_{x\in\mathbb R^2}u(x)>C\,.\end{equation}
Also, we prove in \cite[Lemma~4.4]{kach3} that, for $u\in\mathcal C$
a solution of (\ref{CanonEq}),
\begin{equation}\label{kach3-4.11}
\lim_{x_2\to\pm\infty}\left(\sup_{x_1\in\mathbb
  R}(1-u(x_1,x_2))\right)=0\,.
\end{equation}
Now, let $u_1,u_2\in\mathcal C$ be solutions of (\ref{CanonEq}). We
introduce
\begin{equation}\label{kach3-H-lambda}
\lambda_*=\sup\{\lambda\in[0,1[~:~u_2(x)> \lambda u_1(x)\}\,.
\end{equation}
Then, by (\ref{kach3-4.9}), $\lambda_*>0$. We claim that
$\lambda_*=1$. Once this is shown to hold,
Theorem~\ref{CanonEq-thm} is proved.\\
We argue by contradiction: If $\lambda_*<1$, then
\begin{equation}\label{kach3-w}
\inf_{x\in\mathbb R^2} w(x)=0\,,
\end{equation}
where $w(x)=u_2(x)-\lambda_*u_1(x)$. Now, let
$(x_n)=\left((x_n^1,x_n^2)\right)$ be
a minimizing sequence:
$$\lim_{n\to+\infty}w(x_n)=0\,.$$
Since the maximum principle yields that $w(x)>0$ for all $x$, the
sequence $(x_n)$ should be unbounded, hence we assume that
$\lim_{n\to+\infty }|x_n|=+\infty$. Also, by (\ref{kach3-4.11}), $(x_n^2)$
should be bounded, hence we assume that $\lim_{n\to+\infty}x_n^2=b$.\\
Now, the functions $u_j^n(x_1,x_2)=u_j(x_1+x_1^n,x_2)$, $j=1,2$, solve
(\ref{CanonEq}) in $\mathcal C$, and up to extraction of a
subsequence, they  converge locally to functions\break
$\widetilde u_j\in C_{\rm  loc}^2(\mathbb R\times\{\mathbb
R\setminus[-d,d]\})$, $j=1,2$. Now, $\widetilde u_1$, $\widetilde u_2$
solve (\ref{CanonEq}) in
$\mathcal C$, $\widetilde u_2\geq \lambda_* \widetilde u_1$ and
$\widetilde u_2(0,b)=\lambda_* \widetilde u_1(0,b)$. On the other hand, the
strong maximum principle insures that $\widetilde u_2(x)>\lambda_*
\widetilde u_1(x)$ for all $x\in\mathbb R^2$, hence we have a
contradiction.\hfill$\Box$\\

\section{On the energy of a vortex configuration with pinned vortices}
\label{App-B}

In this section we construct a configuration with several vortices
on the circle $\mathbb S_R^1$, and we discuss the difficulty behind
the estimation of its energy. The motivation to construct such a
configuration is that we  expect vortices of a minimizer of
(\ref{V-EGL}) to be pinned   on the circle $\mathbb S_R^1=\partial
D(0,R)$, and to be uniformly distributed
along this circle.\\
We cover the unit disc $\Omega=D(0,1)$
uniformly by sectors, each having a small opening. Then we define
a periodic induced magnetic field $h$ on the sectors. This
magnetic field will be the natural one
corresponding to a vortex configuration with a single vortex in each
sector. We mention
also that similar constructions have been also introduced
in the papers \cite{AfAlBr, AlBr06, AlBr}.\\
\indent Let $(r,\theta)$ be the polar  coordinates. Let us take
$n=n(\varepsilon)$ equidistant points $(a_i)_{i=1}^N$ on the circle
$\mathbb S^1_R$ with $\theta(a_1)=0$. Upon defining the angular
sector
$$C_1=\left\{(r,\theta)~:~r\in[0,1[,\quad
|\theta|<\frac{2\pi}{n(\varepsilon)}\right\},$$  we get a
decomposition of $D(0,1)$ by a family of disjoint sectors $(C_i)$,
where each $C_i$ corresponds to the point $a_i$ and is congruent to
$C_1$.
\subsection*{\it The test configuration}\ \\
Now we define a measure $\mu$ by:
$$\mu(x)=\left\{\begin{array}{lll}
0&{\rm if}&x\not\in\cup_{i}B(a_i,\varepsilon)\\
\displaystyle\frac2{\varepsilon^2}&{\rm if}&x\in\overline{\cup_i
B(a_i,\varepsilon)},
\end{array}\right.$$
and a function $h'$ in $\Omega=D(0,1)$ by
\begin{equation}\label{V-Def-f'}
\left\{
\begin{array}{rl}
-{\rm div}\,\left(\displaystyle\frac1{u_\varepsilon^2}\nabla h'\right)
+h'=\mu&{\rm in}~\Omega,\\
h'=0&{\rm on}~\partial \Omega.\end{array}\right.\end{equation} We
notice that
$$\int_{C_i}\mu\,\md x=2\pi,\quad\forall~i=1,2,\cdots N\,.$$
We define  an induced magnetic field $h=h'+h_\varepsilon$  (here
$h_\varepsilon$ has been introduced in (\ref{V-hepsilon'})). Then we
define an induced magnetic potential
$A=A'+\frac{1}{u_\varepsilon^2}\nabla^\bot h_\varepsilon$ by
  taking simply
$${\rm curl}\,A'=h'.$$
This choice is always possible as one can take $A'=\nabla^\bot g$
with $g\in H^2(\Omega)$ such that
$\Delta g=h'$.\\
We turn now to define an order parameter $\psi$ which we take in the
form
\begin{equation}\label{V-vsol-psi'}
\psi=u\,u_\varepsilon=\rho\,e^{i\phi}\,u_\varepsilon,\end{equation}
where $\rho$ is defined by:
\begin{equation}\label{V-vsol-rho'}
\rho(x)=\left\{
\begin{array}{cll}
0&{\rm if}&x\in \cup_iB(a_i,\varepsilon),\\
1&{\rm if}&x\not\in\cup_iB(a_i,2\varepsilon),\\
\displaystyle\frac{|x-a_i|}\varepsilon-1&{\rm if}&\exists\,i~{\rm
s.t.}~x\in B(a_i,2\varepsilon) \setminus B(a_i,\varepsilon).
\end{array}\right.
\end{equation}
The phase $\phi$ is defined (modulo $2\pi$) by the relation:
\begin{equation}\label{V-vsol-phi'}
\nabla\phi-A'=-\frac{1}{u_\varepsilon^2}\nabla^\bot h'\quad{\rm
in}~ \Omega\setminus \cup_iB(a_i,\varepsilon),
\end{equation}
and we emphasize here that we do not need to define $\phi$ in
regions where $\rho$ vanishes.

\subsection*{\it The energy of the test configuration}\ \\
To estimate the energy of the test configuration $(\psi,A)$
constructed above, we  express $h'$ by means of a Green's
function, i.e. a fundamental solution of the differential operator
$-{\rm
  div}\,\left(\displaystyle\frac1{u_\varepsilon^2(x)}\nabla\right)+1$.
The existence
and the properties of this function, taken from \cite{AfSaSe, Stamp},
are given in the next lemma.

\begin{lem}\label{Green}
For every $y\in\Omega$ and $\varepsilon\in]0,1]$, there exists a
symmetric function
$\overline\Omega\times\overline\Omega\ni(x,y)\mapsto
G_\varepsilon(x,y)\in\mathbb R_+$ such that~:
\begin{equation}\label{equation-Green}
\left\{\begin{array}{rl}
-{\rm div}\,\left(\displaystyle\frac1{u_\varepsilon^2(x)}\nabla_x
G_\varepsilon(x,y)\right)+G_\varepsilon(x,y)=\delta_y(x)&{\rm in}~\mathcal
D'(\Omega),\\
G_\varepsilon(x,y)\big{|}_{x\in\partial\Omega}=0.
\end{array}\right.
\end{equation}
Moreover, $G_\varepsilon$ satisfies the following properties:
\begin{enumerate}
\item There exists  a constant $C_\varepsilon>0$  such
that
$$0\leq G_\varepsilon(x,y)\leq
C_\varepsilon\left(\big{|}\ln|x-y|\,\big{|}+1\right),\quad\forall~(x,y)\in\overline\Omega\times\overline\Omega\setminus\Delta,$$
where $\Delta$ denotes the diagonal in $\mathbb R^2$.
\item For any compact set $K\subset\Omega$, there exist  constants
$ C>0$ and $\varepsilon_0>0$ such that,
$\forall~\varepsilon\in]0,\varepsilon_0]$,
$$\left|G_\varepsilon(x,y)
+\frac{u_\varepsilon^2(x)}{2\pi}\ln|x-y|\,\right|\leq
C\left\|
\frac{|\nabla u_\varepsilon(x)|}{u_\varepsilon^2(x)}\right\|_{L^\infty(\Omega)}
,\quad\forall~y\in K,~\forall~x\in\overline\Omega.$$
\end{enumerate}
\end{lem}

The field $h'$ can be expressed by means of the function
$G_\varepsilon$,
\begin{equation}\label{h'-Green}
h'(x)=\int_\Omega G_\varepsilon(x,y)\,\mu(y)\,\md
y\,,\quad\forall~x\in\Omega.
\end{equation}
By this expression of the field, one is able in the former
literature to account for the energy resulting from the interaction
of different vortices.\\
By applying Lemma~\ref{Beth-Riv}, one essentially needs to estimate
\begin{equation}\label{Appendix-lem}
\int_\Omega\left(\frac{1}{u_\varepsilon^2(x)}|\nabla
h'|^2+|h'|^2\right)\,\md x=
\int_{\Omega\times\Omega}G_\varepsilon(x,y)\,\mu(x)\mu(y)\,\md x\,\md y.\end{equation}
Unfortunately, due to the rapid oscillation of $u_\varepsilon$,
the decay  of $G_\varepsilon$  in
Lemma~\ref{Green} is not useful to simulate the energy, since the term
$\|\frac{|\nabla u_\varepsilon|}{u_\varepsilon^2}\|_{L^\infty}$ is typically of  order
$\varepsilon^{-1}$.

\section{Remarks on the Ginzburg-Landau functional}
We recall in this appendix some facts taken from \cite{Sand, SaSe}
concerning the Ginzburg-Landau functional~:
\begin{equation}\label{GL-classical}
\mathcal J(u,A)=\int_\Omega\left(
|(\nabla-iA)u|^2+\frac1{2\varepsilon^2}(1-|u|^2)^2+|{\rm
  curl}\,A-H|^2\right)\,\md x,
\end{equation}
defined for functions $(u,A)\in H^1(\Omega;\mathbb C)\times
H^1(\Omega;\mathbb R^2)$. Here $\Omega$ is a two-dimensional simply
connected domain.\\

Let us start by the following result concerning an estimate of the
radius of the set where the modulus of a  function $u$ is small,
provided that $u$ satisfies some energy bound.

\begin{prop}\label{Sandier1}
There exist constants $\alpha_0(\Omega)>0$ and $C>0$ such that, for
any $M>0$, $\delta>0$, $\varepsilon\in]0,1[$ satisfying $\varepsilon
M/\delta^2<\alpha_0$, any $u\in C^2(\overline\Omega;\mathbb C)$ and
$A\in C^0(\Omega;\mathbb R^2)$ satisfying
$$\int_\Omega\left(|\nabla|u|\,|^2
+\frac1{2\varepsilon^2}(1-|u|^2)\right)\,\md x\leq M,$$ one has the
estimate
$$r\big{(}\{x\in\overline\Omega~:~|u(x)|\leq 1-\delta\}\big{)}\leq
C\frac{\varepsilon M}{\delta^2}.$$
\end{prop}

The next proposition provides a lower bound of the functional
$\mathcal J(u,A)$ by means of the degree of $u$ on perforated sets.

\begin{prop}\label{Sandier2}
Let $w$ be a compact subset of $\Omega$. Then for any
$1>\eta>\alpha> r(w)$, there exists a collection of disjoint open
balls $(B(a_i,r_i))_i$ such that
\begin{enumerate}
\item
$\sum_ir_i\leq \eta$\,;
\item
$w\subset \bigcup _i B(a_i,r_i)$\,;
\item For any $H^1$-functions
$u~:\Omega\setminus w\to \mathbb S^1$ and $A~:\Omega \to\mathbb
R^2$,
 we have
$$\int_{B(a_i,r_i)\setminus w}|(\nabla-iA)u|^2\,\md x+
r_i\int_{B(a_i,r_i)}|{\rm curl}\,A-H|^2\,\md x\geq 2\pi|d_i|
\left(\ln\frac{\eta}{\alpha}-\frac{r_i}2\right),$$ where $d_i={\rm
deg}(u,\partial B(a_i,r_i))$ if $B_i\subset \Omega$ and $d_i=0$
otherwise.
\end{enumerate}
\end{prop}

Finally we state a result concerning `Jacobians'.

\begin{prop}\label{Sandier3}
Assume that $M>0$, $R>0$, $\gamma\in]0,1]$ and $\delta\in]0,1/2]$.
Let $u\in C^1(\overline\Omega;\mathbb C)$, $A\in
C^0(\overline\Omega;\mathbb R^2)$ and $(B(a_i,r_i))_i$ a collection
of balls such that
\begin{enumerate}
\item
$|u|\leq 1$ in $\Omega$\,;
\item
$\mathcal J(u,A)\leq M$\,;
\item
$\{x\in\overline\Omega~:~|u(x)|\leq 1-\delta\}\subset \bigcup_i
  B(a_i,r_i)$\,;
\item $\displaystyle\sum_ir_i<R$\,.
\end{enumerate}
Then, there exist measures  $\alpha\in H^{-1}(\Omega)$ and $\beta\in
(C_0^\gamma(\Omega))'$ such that
$$\|\alpha\|_{H^{-1}_0}\leq C\,M\delta^2,\quad\|\beta\|_{(C_0^\gamma)'}
\leq C\,M\,R^\gamma,$$ and, upon putting $d_i={\rm deg}(u,\partial
B(a_i,r_i))$ if $B(a_i,r_i)\subset \Omega$ and $d_i=0$ otherwise, we
have,
$$4\pi\sum_id_i\,\delta_{a_i}-{\rm
  curl}\left[(iu,\nabla_Au)+A\right]=\alpha+\beta.$$
Here $C>0$ is a constant independent of $\varepsilon$, $\delta$, $R$ and $M$.
\end{prop}


\end{document}